\documentclass[10pt]{article}
\usepackage{amsfonts}
\usepackage{amsmath, amsthm, amsfonts, amssymb,mathrsfs,color}
\usepackage{authblk}  
\usepackage{mathtools}
\usepackage[scaled=0.99]{helvet}               
      
\usepackage{orcidlink}
\usepackage{tensor}
\usepackage{hyperref}
\hypersetup{
colorlinks,
citecolor=blue,
filecolor=black,
linkcolor=blue,
urlcolor=blue
}

\usepackage{url}
\usepackage{stix} 
\usepackage{todonotes}
\usepackage{capt-of}

\usepackage{soul}
\usepackage{enumitem}
\usepackage[normalem]{ulem}
\usepackage[numbers,sort&compress]{natbib}

\usepackage
[
a4paper,
vmargin=2.2cm,
hmargin=1.7cm
]{geometry}

\usepackage{indentfirst}

\usepackage{booktabs}
\usepackage{array}
\usepackage{multirow}
\usepackage{makecell}

\def\eps{\varepsilon}

\def\pd{\partial}
\def\F{\mathcal F}
\def\Ex{\mathbb{E}}
\def\d{\,\mathrm{d}}
\def\I{{\mathbf I}}
\def\R{\mathbb{R}}
\def\T{\mathbb{T}}
\def\K{\mathbb K}
\def\W{\mathcal W}
\def\N{\mathbb{N}}
\def\Z{\mathbb{Z}}
\def\D{\mathcal D}
\def\p{\mathbb P}
\def\P{\mathcal P}
\def\B{\mathcal B}
\def\A{\mathcal A}
\def\Q{\mathcal Q}
\def\OP{{\rm OP}}
\def\O{\mathscr O}
\def\S{\mathcal S}

\def\Wlip{W^{1,\infty}}

\def\U{\mathbb U}
\def\LL{\mathscr L}
\def\pas{{\mathbb P}\text{-}{\rm a.s.}}

\def\ol{\overline}
\def\hh{\widehat}

\newcommand{\abs}[1]{\left | #1 \right |}
\newcommand{\norm}[1]{\left\|#1\right\|}
\newcommand{\bk}[1]{ \left(  #1 \right)}

\newcommand{\IP}[1]{\left\langle#1\right\rangle }
\newcommand{\bIP}[1]{\big\langle#1\big\rangle }

\theoremstyle{plain}
\newtheorem{theorem}{Theorem}[section]
\newtheorem{lemma}[theorem]{Lemma}

\theoremstyle{definition}
\newtheorem{definition}{Definition}[section]
\newtheorem{remark}{Remark}[section]

\newtheorem{example}{Example}[section]

\numberwithin{equation}{section}

\newtheorem*{MGoal*}{Main Goal}
\newtheorem*{Explanations*}{Explanations}

\newtheorem{Assumption}{Hypothesis}
\newenvironment{ManualHypo}[1]{
\renewcommand\theAssumption{#1}
\Assumption
}{\endAssumption}

\usepackage{fancyhdr}
\title{{\bf Damping-diffusion-noise interactions in the stochastic Camassa-Holm equation}}
 
\author{
	{\textbf{Diego Alonso-Or\'an}}\thanks{Email: \href{mailto:dalonsoo@ull.edu.es}{dalonsoo@ull.edu.es}}
}
\affil{{\small Departamento de Análisis Matemático y Instituto de Matemáticas y Aplicaciones (IMAULL), Universidad de La Laguna, C/Astrofísico Francisco, Sánchez s/n, 38271 La Laguna, Spain}}

\author{
	{\textbf{Peter H.C.~Pang}}\thanks{Email: \href{mailto:Peter.Pang@nottingham.edu.cn}{Peter.Pang@nottingham.edu.cn}}
}
\affil{{\small School of Mathematical Sciences, University of Nottingham Ningbo China, Ningbo 315000, China}}

\author{
	{\textbf{Hao Tang}}
	\thanks{Email: \href{mailto:haotang@tju.edu.cn}{haotang@tju.edu.cn} (Corresponding author)}
}
\affil{{\small Center for Applied Mathematics, Tianjin University, Tianjin 300072, China}}

\begin{document}

\maketitle

\begin{abstract}
In this paper we investigate the effects of the interaction between time-inhomogeneous damping, non-local diffusion, and noise on classical solutions to the Camassa-Holm equation incorporating these features.  First, a local-in-time theory is established, covering the existence, uniqueness, and a blow-up criterion under relatively general conditions for these interacting mechanisms. Subsequently, we identify different conditions on the interactions between damping, diffusion, and noise that guarantee the global regularity and long-time behaviour of classical solutions. Notably, we demonstrate the existence of an evolution system of measures, which generalizes the concept of invariant measures to the time-inhomogeneous setting.

\end{abstract}


\vskip 6pt


\setcounter{tocdepth}{2}  
\tableofcontents 

\section{Introduction}

The Camassa--Holm (CH) equation  
\begin{equation} 
\partial_t u - \partial_{xxt}u+ 3u \partial_xu = 2\partial_xu\partial_{xx}u + u\partial_{xxx}u,\ \ 
u=u(t,x),\label{CH} 
\end{equation} 
was first introduced in the context of hereditary symmetries studied by Fuchssteiner and Fokas \cite{Fuchssteiner-Fokas-1981-PhyD} as a bi-Hamiltonian generalization of the Korteweg-de Vries equation. 
Camassa and Holm later rediscovered it 
\cite{Camassa-Holm-1993-PRL} as a model for the unidirectional 
propagation of shallow water waves over a flat 
bottom. For a rigorous justification of the derivation of equation \eqref{CH}, see \cite{Constantin-Lannes-2009-ARMA}.

Extensive studies have been conducted on the CH equation \eqref{CH}, and it is beyond the 
scope of this paper to provide a comprehensive review. We briefly mention a few related results. The Cauchy problem, which addresses both local and global well-posedness for the CH equation, along with various solution concepts—such as strong, weak, conservative, and dissipative—has been thoroughly investigated. For example, \eqref{CH} is locally well-posed for strong solutions in the Sobolev space $H^s$ for $s > 3/2$ (see \cite{Li-Olver-2000-JDE, Constantin-Escher-1998-CPAM}). The phenomenon of solution blow-up in \eqref{CH}, characterized by wave-breaking—where the solution remains bounded while its slope becomes unbounded in finite time—has been extensively explored (see \cite{Constantin-2000-JNLS, Constantin-Escher-1998-CPAM, McKean-1998-AJM}). 
Moreover, the CH equation allows for global weak solutions within the $H^1$ regularity class, where the solutions can be either conservative or dissipative, depending on the method of extension beyond the blow-up time. For further details, see \cite{Bressan-Constantin-2007-AA, Bressan-Constantin-2007-ARMA, Holden-Raynaud-2007-CPDE, Holden-Raynaud-2009-DCDS} and references therein.

Let $\I$ be the identity operator. \eqref{CH} can be 
reformulated into
\begin{equation*}
\pd_t u + u\,\pd_x u + F(u)=0
\end{equation*}
with 
\begin{equation}\label{F define}
F(u)\triangleq\pd_x \bk{\I -\pd^2_{x}}^{-1} 
\bk{u^2 + \frac12 \bk{\pd_x u}^2}.
\end{equation}
When  $x\in\R$ or $x\in\T\triangleq\R/2\pi\Z$,  the operator 
$(\I -\pd^2_{x})^{-1}$ can be 
defined by Fourier transform (see \eqref{Ds define} below).

\subsection{Motivations and main goal}

As noted previously, solutions to \eqref{CH} can blow up in the form of wave-breaking within finite time. It is indeed feasible to precisely predict the onset of wave-breaking by establishing a necessary and sufficient condition for solutions to the Cauchy problem associated with \eqref{CH} (see \cite{McKean-1998-AJM}). Consequently, various mechanisms have been proposed in the literature to  \textit{regularise}  the CH equation.

From a deterministic perspective, regularisation is typically achieved by introducing viscosity or damping terms into the equation. Viscosity acts as a dissipation mechanism that smooths the solution and enhances stability, while damping reduces the energy of the system, controlling the growth of perturbations.  For CH-type equations, viscous  regularisations have been studied in one- and multi-dimensional settings, see \cite{Foias-Holm-Titi-2002-JDDE, Coclite-Holden-Karlsen-2005-DCDS}. 
The damped CH-type equations have also been explored in \cite{Lenells-Wunsch-2013-JDE, Wu-Yin-2008-SIAM, Wu-Yin-2009-JDE}. These models can be unified as follows:
\begin{align*} 
\pd_t u + \eps (-\pd_{x}^2)^{\theta}u 
+\lambda(t) \, u  +u \,\pd_xu +F(u) =0,
\end{align*}
where $\eps (-\pd_{x}^2)^{\theta}u$ represents the (non-local) diffusion term with $\theta \in (0,1]$, $\lambda(t) \, u$ is the damping term, and $F(u)$ is defined in \eqref{F define}. When $\lambda(t) \equiv \lambda_0 > 0$, the term $\lambda_0 u$ is commonly referred to as the \textit{weakly dissipative} term, since it induces the term $\lambda_0 (u - \pd^2_{x} u)$, which is added to \eqref{CH}.

From a stochastic perspective, various stochastic versions of \eqref{CH} have been extensively studied in recent years, with regularisation effects due to random perturbations being  recognised. In the case of (It\^o) multiplicative noise, linear noise can ensure that a strong solution with small initial data exists globally with high probability (see \cite{Rohde-Tang-2021-JDDE, Chen-Miao-Shi-2023-JMP}), while certain faster-growing non-linear noise can almost surely prevent blow-up (see \cite{Rohde-Tang-2021-NoDEA, Tang-Yang-2023-AIHP, Tang-2023-JFA}). For additional studies on the stochastic CH-type equation with (It\^o) multiplicative noise, refer to \cite{Tang-2018-SIMA, Chen-Duan-Gao-2021-CMS, Ren-Tang-Wang-2024-POTA, Chen-Duan-Gao-2023-AAP, Chen-Miao-Shi-2023-JMP}. When the noise is of (Stratonovich) transport-type, we refer to \cite{Albeverio-etal-2021-JDE, Alonso-Rohde-Tang-2021-JNLS, Holden-Karlsen-Pang-2023-DCDS, Galimberti-etal-2024-JDE}.

However, to the best of the authors' knowledge, there is \textbf{no} existing work that studies the interaction among these three phenomena in the stochastic CH equation. Therefore, the main goal of this paper is
\begin{MGoal*}
\textit{To compare the effects of non-local diffusion, damping, and noise on classical solutions}.
\end{MGoal*}

To this end, 
we fix the noise structure considered in advance.   Let $\{W_k,\widetilde W_k\}_{k\ge1}$ be mutually independent
real-valued two-sided Brownian motions on a complete probability
space $(\Omega, \{\mathcal F_t\}_{t\in\mathbb R},\mathbb P)$, where
$W_k(0)=\widetilde W_k(0)=0$ and
$\{\mathcal F_t\}_{t\in\mathbb R}$ is the usual augmentation of
the filtration generated by their increments up to time $t$. $\circ {\d}W_k(t)$ denotes the Stratonovich stochastic differential, and ${\rm d}\widetilde{W}_k(t)$ is the It\^o stochastic differential.  In this paper, we examine the following noise structure:
$$\sum_{k = 1}^\infty \Q_k u(t) \circ{\rm d} W_k(t) +\sum_{k = 1}^\infty h_k(t,u(t))\d \widetilde{W}_k(t),$$
where  $h_k(t,u)$ is a sequence of non-linear functions and $\{\Q_k\}_{k\ge1}$ is a sequence of pseudo-differential operators.

We note that such a noise structure encompasses the (Stratonovich) transport-type noise $\Q_k u\circ \mathrm{d} W_k$ (if $\Q_k$ reduces to the classical derivative operator of order 1) and the (It\^o) multiplicative noise $h_k(t,u)\d \widetilde{W}_k$.

Let $F(u)$ be given in \eqref{F define}. In this paper, we will consider the following stochastic CH equation on $\K$, with a focus on the interactions between non-local diffusion, damping, and noise:
\begin{align}\label{SCH}
\d u(t) + \Big\{\eps \Lambda^{2\theta}u +\lambda(t) \, u  +u \,\pd_xu +F(u)\Big\}(t)\d t  
= \sum_{k = 1}^\infty \Q_k u(t) \circ{\rm d} W_k(t) + \sum_{k = 1}^\infty h_k(t,u(t))\d \widetilde{W}_k(t),
\end{align}
where 
\begin{itemize} 
\setlength\itemsep{0em}
\item $\eps\Lambda^{2\theta}$, with $\Lambda\triangleq\bk{-\pd^2_{x}}^{1/2}$, $\eps\ge0$, and $\theta\in(0,1]$, is a non-local diffusive operator (see \eqref{Lambda s define} below); 
\item $\lambda(t)\, u$ is the damping term with a time-dependent coefficient $\lambda\in L^1_{\rm loc}(\R; [0,\infty))$;
\item $\{\mathcal{Q}_k\}_{k\ge1}$ is a sequence of pseudo-differential operators (see Section \ref{Section: Background} below).
\end{itemize}

We now briefly discuss the \textit{noise structure} $\sum_{k = 1}^\infty \Q_k u(t) \circ{\rm d} W_k(t)$, which generalises transport-type noise in a non-local way. In \cite{Mikulevicius-Rozovskii-2004-SIAM, Mikulevicius-Rozovskii-2005-AoP}, randomness was introduced into the system of equations at the Lagrangian level by imposing a stochastic forcing $\sigma \circ {\,\mathrm{{d}}}W(t)$ (where $W(t)$ is a standard 1-D Brownian motion) in the equation for the characteristic lines $X = X(t,x)$ with an undetermined local velocity $u$. Specifically, this was done via the Lagrangian stochastic differential equation:
\begin{equation*}
{\rm d}X(t,x) = u(t, X(t,x)) {\,\mathrm{{d}}}t + \sigma (t,X(t,x)) \circ {\mathrm{{d}}}W(t),\quad X(0,x)=x.
\end{equation*}
By reformulating the system in terms of $u(t,X(t,x))$, we obtain what is nowadays known as transport noise (cf. \cite{Holm-2015-ProcA}). In the one-dimensional case, transport noise refers to noise of the following structure:
\begin{equation*}
\sum_{k = 1}^\infty \Q_k u(t) \circ{\rm d} W_k(t),\quad \Q_k = f_k \partial_x + g_k \I,\quad
f_k\ \text{and}\ g_k\ \text{are sufficiently regular functions}.
\end{equation*}
After the works \cite{Mikulevicius-Rozovskii-2004-SIAM, Mikulevicius-Rozovskii-2005-AoP, Holm-2015-ProcA}, many subsequent investigations studying the properties of equations in fluid dynamics with transport-type noise have appeared in the literature (cf. \cite{Alonso-Bethencourt-2020-JNLS, Alonso-etal-2019-NODEA, Alonso-Rohde-Tang-2021-JNLS, Albeverio-etal-2021-JDE, Crisan-Holm-2018-PhyD, Crisan-Flandoli-Holm-2019-JNLS, Holden-Karlsen-Pang-2021-JDE,Lang-Crisan-2022-SPDE, Holden-Karlsen-Pang-2023-DCDS}).

In contrast to these previous works, in this paper, we consider the following case:
\begin{equation*}
\mathcal{Q}_{k}\ \text{is a pseudo-differential operator},
\end{equation*}
which significantly extends previous results. We refer to noise of the structure $\mathcal{Q}_{k}u \circ {\,\mathrm{d}} W_{k}$ as \textit{pseudo-differential noise} (see Section \ref{Section: Background} for the definition and background on pseudo-differential operators). Such an extension to pseudo-differential noise $\mathcal{Q}_{k}u \circ {\,\mathrm{d}} W_{k}$ also enables the modelling of non-local random interactions. Since pseudo-differential operators extend classical differential operators in a non-local manner, studying pseudo-differential noise provides a more flexible framework for modelling complex phenomena that involve long-range random interactions. This can be particularly useful in turbulence models, where the behaviour of fluid at one point is influenced by the behaviour of fluid at distant points (for example, see \cite{Bernard-Erinin-2018-JFM} for non-local Reynolds shear stress and \cite{Hamba-2022-JFM} for non-local eddy diffusivity). To the best of the authors' knowledge, general pseudo-differential noise in the Stratonovich sense has only been studied very recently in \cite{Tang-Wang-2022-arXiv, Tang-2023-JFA}.

Before summarizing the main results of this work, it is convenient to rewrite  \eqref{SCH}. Recalling that for a continuous semi-martingale $g(t)$, we have 
$g(t)\circ \d W_k(t)= g(t)  \d W_k(t)+ \frac{1}{2}  
\d \left\langle g, W_k\right\rangle (t),$
$k\ge1$,
where $\left\langle\cdot,\cdot\right\rangle$ stands for the quadratic variation. We are thereby motivated to study the following Cauchy problem of \eqref{SCH}:
\begin{equation}\label{Cauchy problem-Ito}
\left\{\begin{aligned}
\d u(t)&+\Big\{\eps\Lambda^{2\theta} u 
+\lambda(t)\, u+u\,\pd_x u +F(u)\Big\}(t) \d t\\
&=\frac12\sum_{k=1}^{\infty}{\Q}^2_ku(t) \d t  
+\sum_{k=1}^{\infty} {\Q}_ku(t) \d  W_k(t)
+\sum_{k=1}^{\infty}h_k(t,u(t))\d \widetilde{W}_k(t),
\quad t>0,\\
u(0)&= u_0.
\end{aligned} 
\right.
\end{equation}

\subsection{Summary of the main results}
\label{sec:mainres_summary}

For the reader's convenience, the main results of this paper are briefly outlined in the following table, which provides an overview of the different scenarios considered and their corresponding outcomes. We emphasize that the table serves only as a schematic guide; the precise statements of the results are provided in Section \ref{sec:mainres_precise}. A  detailed explanation of the key difficulties encountered and the novel approaches developed to overcome them is provided in Section \ref{Sect:remarks}.
Let $\K=\R$ or $\K=\T=\R/2\pi\Z$. 

\begin{table}[h]
\centering
\begin{tabular}{
>{\raggedright\arraybackslash}p{4.7cm}
>{\raggedright\arraybackslash}p{10.3cm}
}
\toprule
\textsc{\textbf{Results}} & \textsc{\textbf{Conditions on Damping-Diffusion-Noise Interactions}} \\
\midrule

\textbf{Local-in-Time Theory on $\K$ (Theorem \ref{Thm:local-theory}):} 
\begin{itemize}[leftmargin=*] 
\setlength\itemsep{3.5em}
\item[]\textit{Existence \& Uniqueness \& Blow-up Criterion}
\end{itemize}
& 
\vspace*{0.8cm}
\begin{itemize}[leftmargin=*,nosep]

\item Diffusion may be \textit{degenerate} or \textit{weak} ($\varepsilon \geq 0$, $\theta \in (0,1]$)
\item $\sum_{k=1}^\infty h_k(t,\cdot)$ is locally Lipschitz with growth condition permitting \textit{high nonlinearity}
\item Damping may be \textit{degenerate} ($\lambda(t) \geq 0$)
\end{itemize} \\

\midrule

\textbf{Global-in-Time Theory on $\K$ (Case 1, Theorem \ref{Thm: global 1}):}
\begin{itemize}[leftmargin=*] 
\setlength\itemsep{3.5em}
\item[] \textit{Global existence}

\item[] \textit{Solution decay}
\end{itemize}

& 
\vspace*{0.8cm}
\begin{itemize}[leftmargin=*,nosep]
\item Diffusion satisfies $\varepsilon > 0$, $\theta \in (1/2,1]$
\item $\sum_{k=1}^\infty h_k(t,\cdot)$ has \textit{linear growth}
\item Damping may be \textit{degenerate} ($\lambda(t) \geq 0$)
\vspace*{0.55cm}
\item Additional requirement: Damping $\lambda(t)u$ dominates noise effects from $\sum_{k=1}^\infty (\Q_k u \circ \mathrm{d} W_k + h_k(t,u) \mathrm{d} \widetilde{W}_k)$
\end{itemize} \\

\midrule

\textbf{Global-in-Time Theory (Case 2, Theorem \ref{Thm: global 2}):}
\begin{itemize}[leftmargin=*] 
\setlength\itemsep{3.5em}
\item[] \textit{Global existence on $\K$}

\item[] Existence of evolution system of measures on $\T$
\end{itemize}

& \vspace*{0.8cm}
\begin{itemize}[leftmargin=*,nosep]
\item $\sum_{k=1}^\infty h_k(t,\cdot)$ has \textit{sufficiently fast growth}
\item Diffusion may be \textit{degenerate} or \textit{weak} ($\varepsilon \geq 0$, $\theta \in (0,1]$)
\item Damping may be \textit{degenerate} ($\lambda(t) \geq 0$)
\vspace*{0.55cm}
\item Additional requirement: Appropriate interaction between noise terms $\sum_{k=1}^\infty (\Q_k u \circ \mathrm{d} W_k + h_k(t,u) \mathrm{d} \widetilde{W}_k)$ and damping $\lambda(t)u$
\end{itemize} \\

\bottomrule
\end{tabular}
\end{table}

\textbf{Plan of the paper.} 
We conclude this introduction by outlining the structure of this paper. Section~\ref{Section: Background} introduces the necessary notation and background on pseudo-differential operators. Section~\ref{Section:definition hypotheses} presents the definitions of solutions and the evolution system of measures, along with the core hypotheses used in this work. Section~\ref{sec:mainres_precise} provides the precise statements of our main results (Theorems~\ref{Thm:local-theory}--\ref{Thm: global 2}), which were summarized earlier in Section~\ref{sec:mainres_summary}. The main difficulties and the novelties of our results/approaches, compared with earlier research, are discussed in Section~\ref{Sect:remarks}. Section~\ref{sec:aux_lemmas} compiles several technical estimates essential for the subsequent proofs. The proofs for Theorems~\ref{Thm:local-theory}--\ref{Thm: global 2} are detailed in Sections~\ref{sec:localtheory}--\ref{Section:Diffusion-damping-noise-2}: Section~\ref{sec:localtheory} establishes the local-in-time theory, Section~\ref{Section:Diffusion-damping-noise-1} addresses the first global existence result concerning the diffusion-damping-noise interplay, and Section~\ref{Section:Diffusion-damping-noise-2} proves global existence under fast growing noise and the existence of an evolution system of measures.

\section{Background, main results and remarks}

\subsection{Notations and background}\label{Section: Background}

We begin by laying out the notation
used subsequently. 
Unless otherwise stated (as in Section 
\ref{Section:T2:invariant measures}), 
the analysis of \eqref{SCH} applies 
to both the whole space and periodic cases. As before, 
 we denote $$\K=\R\ \ \text{or}\ \ \K=\T=\R/2\pi\Z$$
to unify the analysis. We use $\lesssim$ 
and $\gtrsim$ to indicate that estimates hold 
up to universal \textit{deterministic} constants. 
For linear operators $\mathcal{A}$ and $\mathcal{B}$, 
$[\mathcal{A},\mathcal{B}]
\triangleq\mathcal{A}\mathcal{B}
-\mathcal{B}\mathcal{A}$.  The $L^{2}$-adjoint 
operator of $\mathcal{A}$ is denoted by $\mathcal{A}^{*}$.

For two topological spaces $E_1$ and $E_2$, 
$\mathscr{B}(E_1;E_2)$ is  the class of all 
measurable maps from $E_1$ to $E_2$ and 
$C(E_1;E_2)$ denotes the class of all
continuous maps from $E_1$ to $E_2$. 
$\mathscr{B}_{B}(E_1;E_2)$ and   $C_{B}(E_1;E_2)$ are the set of bounded elements   in $\mathscr{B}(E_1;E_2)$ and $C(E_1;E_2)$, respectively.   Besides, $C_{B,U}(E_1;E_2)$ is the set of bounded uniformly continuous functions  from $E_1$ to $E_2$.  
When $E_2=\R$, we simply write  
$C(E_1)=C(E_1;\R)$ etc.
When $E_1$ and $E_2$ are separable Hilbert spaces, 
$\mathscr L(E_1;E_2)$ is the class of bounded 
linear operators from $E_1$ to $E_2$ and
$\mathscr L_{2}(E_1;E_2)$ is the class of 
Hilbert-Schmidt operators between them.
For $-\infty<a<b<\infty$, we define
$$ \intbar_a^b f(t)\d t\triangleq\frac{1}{b-a}\int_{a}^{b}f(t)\d t.$$

When stating the assumptions on the noise, 
we will use functions in the following classes:
\begin{equation}\label{SCRK} 
\mathscr K\triangleq\left\{
\begin{aligned}
&K(x,y)\in\mathscr{B}
\big( \mathbb R\times [0,\infty);
(0,\infty)\big):\\
K(x,y) \
&\text{is locally integrable in}\ x
\text{\ and increasing in}\ y
\end{aligned}
\right\},
\end{equation} 
and 
\begin{equation}\label{SCRV}
\mathscr{V}\triangleq\left\{f\in C^2([0,\infty);[1,\infty)):
f'(x)>0,\   f''(x)\le 0,\ \lim_{x\to\infty} f(x)=\infty\right\}.
\end{equation}

When $E$ is a topological space,  
$\mathfrak B(E)$ denotes the Borel sets of $E$ and 
$\mathbf P(E)$ is the collection of Borel 
probability measures on $(E,\mathfrak B(E))$.

The Fourier transform $\mathscr F_{x\to\xi}$ 
and inverse Fourier transform 
$\mathscr F^{-1}_{\xi\to x}$ on $\R$ are defined by
$$
(\mathscr F_{x\to\xi}f)(\xi)
\triangleq\int_{\mathbb R}f(x){\rm e}^{-{\rm i}x\xi}\d x,\ \ 
(\mathscr F_{\xi\to x}^{-1}f)(x)
\triangleq  \frac{1}{2\pi}\int_{\mathbb R}
f(\xi){\rm e}^{{\rm i} x\xi}\d \xi.
$$ 
On torus, i.e., $x\in\T$, the Fourier transform 
$\mathscr F_{x\to k}$ and inverse Fourier 
transform $\mathscr F^{-1}_{k\to x}$ are given by
$$
(\mathscr F_{x\to k}f) (k)
\triangleq\int_{\mathbb T}f(x){\rm e}^{-{\rm i}xk}\d x,\ \ 
(\mathscr F_{k\to x}^{-1}f)(x)
\triangleq  \frac{1}{2\pi}\sum_{k\in\Z}f(k){\rm e}^{{\rm i} xk}.
$$ 
For any $s\in\R$, we define 
\begin{equation}\label{Lambda s define}
\Lambda^s=(-\pd^2_{x})^{s/2}\triangleq\left\{
\begin{aligned}
&\mathscr F^{-1}_{\xi\to x}\Big[|\xi|^{s} 
\big(\mathscr F_{x\to\xi}\cdot\big)(\xi)\Big] 
\ \text{on}\ \R\\
&\mathscr F^{-1}_{k\to x}\Big[|k|^{s} 
\big(\mathscr F_{x\to k}\cdot\big)(k)\Big] \  
\text{on}\ \T,\\
\end{aligned}
\right.
\end{equation} 
and
\begin{equation}\label{Ds define}
\D^s=(\I-\pd^2_{x})^{s/2}\triangleq\left\{
\begin{aligned}
&\mathscr F^{-1}_{\xi\to x}\Big[(1+|\xi|^2)^{s/2} 
\big(\mathscr F_{x\to\xi}\cdot\big)(\xi)\Big] 
\ \text{on}\ \R\\
&\mathscr F^{-1}_{k\to x}\Big[(1+|k|^2)^{s/2} 
\big(\mathscr F_{x\to k}\cdot\big)(k)\Big] \  
\text{on}\ \T.\\
\end{aligned}
\right.
\end{equation} 
For $s\ge0$, we define the Sobolev spaces $H^s$ on $\K$ 
with values in $\R$ as 
$$H^s(\mathbb R;\R)
\triangleq\ol{C_0^\infty(\mathbb R;\R)}^{\|\cdot\|_{H^s}},\quad
H^s(\mathbb T;\R)
\triangleq \ol{C^\infty(\mathbb T;\R)}^{\|\cdot\|_{H^s}},
$$
where the norm $\|\cdot\|_{H^s}$ is induced by 
the inner product $\IP{f,g}_{H^s}\triangleq\IP{\D^s f ,\D^s g}_{L^2}.$
Let $W^{1,\infty}(\mathbb K;\R)$ be the set of 
weakly differentiable functions $f:\mathbb K\to\R$ with
$$\|f\|_{W^{1,\infty}}\triangleq 
\sum_{\alpha=0,1}\|\partial_x^{\alpha} f \|_{L^\infty} <\infty.$$
For $s\ge0$ and $p\in [1,\infty]$, we will simply write
$$H^s=H^s(\mathbb K;\R),\ \ 
W^{1,\infty}=W^{1,\infty}(\mathbb K;\R), \ \ 
L^p=L^p(\mathbb K;\R).$$
In particular, we let $H^\infty\triangleq\cap_{s\ge0} H^s$.   
In the following we will  also need to equip $H^s$ with a non-standard topology:
\begin{equation}\label{X-s1s2}
X_{s_1,s_2}\triangleq\big(H^{s_1},\mathscr{T}^{s_2}\big),\quad s_1 \ge s_2,\quad
\mathscr{T}^{s_2}\ \text{is the topology induced by}\ \|\cdot\|_{H^{s_2}}.
\end{equation}
If $s_1>s_2$, the space $X_{s_1, s_2}$ is not complete and 
the completion $\ol{X_{s_1,s_2}}^{\mathscr{T}^{s_2}}$  is $H^{s_2}$.

Let $\mathbb{N}_0 \triangleq \mathbb{N} \cup \{0\}$.
We now introduce the pseudo-differential operators. 
For any $s\in \R$, the symbol class 
$\S^s(\R\times \R)\subset
C^\infty (\R\times \R;\mathbb{C})$ is defined by 
\begin{align*}
 \S^s(\R\times \R;\mathbb{C})\triangleq 
&\left\{\mathfrak{p} \, : \, 
\forall\, \beta,\alpha\in \N_{0},\ \exists\, C(\beta,\alpha)>0\ 
\text{such that}\ \sup_{(x,\xi)\in\R \times \R}
\frac{\left|\partial_{x}^{\beta} \partial_{\xi}^{\alpha} 
\mathfrak{p} (x, \xi)\right|}{(1+|\xi|)^{ s-\alpha}}<C(\beta,\alpha)\right\}.
\end{align*} 
It is well-known that $\S^{s}(\mathbb R\times \mathbb R;\mathbb{C})$ is a Fr\'{e}chet space equipped with the topology generated by the seminorms $\{|\cdot |^{\beta , \alpha ;s}_{\mathbb R\times \mathbb R}\}_{
\beta , \alpha \in \mathbb N_{0}}$, where
\begin{equation*}
|\mathfrak{p}|^{\beta , \alpha ;s}_{\mathbb R\times \mathbb R}\triangleq
\max_{0\le j\le \beta,\, 0\le \ell\le \alpha}
\sup _{(x,\xi )\in \mathbb R\times \mathbb R} {|\partial _{x}^{
j} \partial _{\xi}^{\ell} {\mathfrak{p}}(x, \xi )|}{(1+|
\xi |)^{ -s+\ell}}.
\end{equation*}
For any $\alpha\in \mathbb N_0$, the  partial 
difference operator $\triangle_k^\alpha$ is  
$$
(\triangle_k^\alpha g)(k)  
\triangleq\displaystyle\sum_{\gamma\in\N_0, \gamma\leq \alpha} 
(-1)^{\alpha-\gamma}\binom{\alpha}{\gamma}g(k+\gamma),
\ \ g:\mathbb{Z}\to \mathbb{C},\ \  k\in \mathbb Z.$$
Then the (toroidal) symbol class of order $s$ for 
$s\in\R$ is defined as (cf. \cite{Ruzhansky-Turunen-2010-Book}):
\begin{align*} 
\S^s(\T\times \Z;\mathbb{C})\triangleq 
 \left\{\mathfrak{p} \, : \, 
\begin{aligned}
&\mathfrak{p}(\cdot,k)\in  C^\infty (\T;\mathbb{C}) \ \text{for all}\ k\in\Z;\\
\forall\, \beta,\alpha\in \N_0,\ &\exists\, C(\beta,\alpha)>0\ 
\text{such that}\ \sup_{(x,k)\in \T\times \Z} 
\frac{\left|\partial_{x}^{\beta}\triangle_{k}^{\alpha} 
\mathfrak{p} (x, k)\right|}{(1+|k|)^{ s-\alpha}}<C(\beta,\alpha)
\end{aligned}
\right\}.
\end{align*} 
Similarly,  $\S^{s}(\mathbb T\times \mathbb Z;\mathbb{C})$ is  a Fr\'{e}chet space with the seminorms $\{|\cdot |^{\beta , \alpha ;s}_{\mathbb T\times \mathbb Z}\}_{
\beta , \alpha \in \mathbb N_{0}}$ defined by
\begin{equation*}
|\mathfrak{p}|^{\beta ,\alpha ;s}_{\mathbb T\times \mathbb Z}\triangleq
\max_{0\le j\le \beta,\, 0\le\ell\le \alpha}
\sup _{(x,k)\in \mathbb T\times \mathbb Z} 
{|\partial _{x}^{j}\triangle _{k}^{\ell} \mathfrak{p}(x, k)|}{(1+|k|)^{-s+\ell}}.
\end{equation*}
The pseudo-differential operator 
with symbol $\mathfrak{p}$ is defined by
\begin{equation}
\OP(\mathfrak{p})\triangleq\P,\label{OP define 1}
\end{equation}
\begin{align}
[\P f](x)\triangleq 
\begin{cases}
\left\{\mathscr F^{-1}_{\xi\to x}\Big[\mathfrak{p}(z,\xi) 
\big(\mathscr F_{x\to\xi}f\big)(\xi)\Big]\right\}_{z=x},\ \ 
&\text{if}\ \mathfrak{p}\in 
\S^s(\R\times \R;\mathbb{C}),\vspace*{4pt}\\
\left\{\mathscr F^{-1}_{k\to x}\Big[\mathfrak{p}(z,k) 
\big(\mathscr F_{x\to k}f\big)(k)\Big]\right\}_{z=x},
\ \ &\text{if}\ \mathfrak{p} \in
\S^s(\T\times \Z;\mathbb{C}).
\end{cases}
\label{OP define 2}
\end{align}

Let $\overline{g}$ denote the complex conjugate of $g\in\mathbb C$. 
Throughout this paper, pseudo-differential operators 
are assumed to be real-valued,  i.e., when $f$ is real, 
$[\OP({\mathfrak{p}})f]$ is also real. 
Equivalently, we require  
\begin{align}
\mathfrak{p}(x,-\xi)=\overline{ \mathfrak{p}(x,\xi)}\ \ \text{if}\ \
(x,\xi)\in\R \times \R ,\label{Real symbol R}\\
\mathfrak{p}(x,-k)=\overline{ \mathfrak{p}(x,k)}\ \  \text{if}\ 
\ (x,k)\in\T  \times \Z .\label{Real symbol T}
\end{align}

We denote by $\S^s\big(\T\times \R;\mathbb C\big)$ 
the class of symbols $\mathfrak{p}\in 
\S^s\big(\R\times \R;\mathbb C\big)$ 
such that $\mathfrak{p}(\cdot,\xi)$ is 
$2 \pi$-periodic for all $\xi\in\R$. 
Moreover, according to 
\cite[Corollary 4.5.7]
{Ruzhansky-Turunen-2010-Book},  
any bounded  set  in
$\S^s\big(\T\times \Z;\mathbb C\big)$ 
(in the sense of boundedness in Fr\'{e}chet 
space) coincides 
with the restriction of a bounded  set in 
$\S^s\big(\T\times \R;\mathbb C\big)$. 
Therefore, when context precludes ambiguity, we write 
\begin{equation*}
\S^s\triangleq
\left\{\mathfrak{p}\in \S^s(\R\times \R;\mathbb{C})\,  
: \,   \eqref{Real symbol R}\
\text{holds} \right\}\  \text{or}\ 
\left\{\mathfrak{p}\in\S^s(\T\times \Z;\mathbb{C})\,  
: \,  \eqref{Real symbol T}\ \text{holds} \right\}.
\end{equation*}
Let $\OP(\cdot)$  be given in 
\eqref{OP define 1}--\eqref{OP define 2}.
Via the above-mentioned relation between 
symbol classes $S^{s}\big(\mathbb T
\times \mathbb Z; \mathbb C\big)$ 
and $S\big(\mathbb T
\times \mathbb R;\mathbb C\big)$,  
we also simply write
\begin{align*}
\OP\S^s\triangleq\Big\{\OP({\mathfrak{p} })\, 
: \,  {\mathfrak{p} } \in\S^s\Big\},\quad
|\cdot |^{\beta ,\alpha ;s}\triangleq|\cdot |^{\beta ,\alpha ;s}_{\mathbb R\times \mathbb R}\ \text{or} \ |\cdot |^{\beta ,\alpha ;s}_{\mathbb T\times \mathbb Z},
\quad  
s\in\R
\end{align*}
to unify ${\mathrm{OP}}\mathcal S^{s}
\big(\mathbb R \times \mathbb R;
\mathbb C\big)$ and ${\mathrm{OP}}\mathcal S^{s}
\big(\mathbb T
\times \mathbb Z;\mathbb C\big)$ when ambiguity is unlikely. 
In the following, we will also consider 
symbols only depending on the frequency 
variable $\xi$ (if $x\in\R$) or $k$ (if $x\in\T$). 
To highlight the differences, we let
\begin{equation*}
\S_0^s\triangleq\left\{\mathfrak{p}\in \S^s\, : \,  
\begin{aligned}
&\mathfrak{p}(x,\xi)=\mathfrak{p}(\xi),\ \ \text{if}\ 
(x,\xi)\in\R\times \R\\
&\mathfrak{p}(x,k)=\mathfrak{p}(k),\ \, \text{if}\ 
(x,k)\in\T\times \Z
\end{aligned}
\right\},
\end{equation*}
and similarly define
\begin{align*}
\OP\S_0^s\triangleq \Big\{\OP({\mathfrak{p} })\, : 
\,  {\mathfrak{p} } \in \S_0^s\Big\},\ \ s\in\R.
\end{align*}

By \cite{Abels-2012-Book,Taylor-1991-book}, 
for $q,s\in\R$ and $\mathfrak{q}\in\S^s$,  there
are constants $\widetilde \beta ,\widetilde \alpha \in \mathbb N_{0}$ 
and 
$C=C(s,q)>0$ such that
\begin{align}
\|{\mathrm{OP}}(\mathfrak{q})\|_{\mathscr L(H^{q+s};H^{q})}\leq C(s,q)|
\mathfrak{q}|^{\widetilde \beta ,\widetilde \alpha ; s}.\label{OP boud Hs}
\end{align}
A subset $\O$ of $\OP\S^s$ is called 
bounded if $\big\{\mathfrak{p}: \OP(\mathfrak{p})\in\O\}\subset \S^s$ 
is bounded.  
Therefore, if $\O\subset \OP\S^s$ is  bounded, 
\eqref{OP boud Hs} tells us that for some $C_0(s,q)>0$,
\begin{align*}
\sup_{\mathcal P\in\O}\|\mathcal P f\|_{H^{q}}
\leq C_0(s,q)\|f\|_{H^{q+s}}, \ \ f\in H^{q+s}.
\end{align*}

\subsection{Definitions and hypotheses}\label{Section:definition hypotheses}

In this section we list two important notions: the notion of pathwise classical solution to \eqref{Cauchy problem-Ito} and the notion of an evolution system of measures, and make the precise assumptions on noise coefficients. 

To begin with, we present the following definition of a pathwise classical solution to \eqref{Cauchy problem-Ito}:
\begin{definition}\label{pathwise solution definition}
Let $u_0$ be an $H^s(\K)$-valued 
$\mathcal{F}_{0}$-measurable  random variable with $s>3/2$.  
\begin{enumerate}
\item A local pathwise classical solution to 
\eqref{Cauchy problem-Ito} is a pair $(u,\tau)$, 
where $\tau$ is a stopping time satisfying 
$\p(\tau>0)=1$ and
$(u(t))_{t\in[0,\tau)}$ is 
$\mathcal{F}_t$-progressively measurable such that 
$$
\sup_{t'\in[0,t]}\|u(t')\|_{H^{s}}<\infty,\ t\in[0,\tau)\ \ \pas,
$$ 
and the following 
equation holds in $C^1(\K)$:
\begin{align*} 
u(t)-u_0
\ + &\int_0^{t}
\Big[\eps\Lambda^{2\theta}u + \lambda(t')\, u
+u\,\pd_x u +F(u)-\frac12\sum_{k=1}^{\infty}\Q^2_ku\Big](t')\d t'\\
\ = &\int_0^{t}\sum_{k=1}^{\infty}
\left(\Q_ku(t') \d  W_k(t')+h_k(t',u(t'))\d 
\widetilde{W}_k(t')\right),\quad t\in[0,\tau)\quad \pas
\end{align*}
\item  Additionally,  a local solution $(u,\tau^*)$ 
is called maximal if $\tau^*>0$ almost surely and 
\begin{equation*}
\limsup_{t\to\tau^*}\|u(t)\|_{H^s}=\infty\ \ 
\text{a.s. on}\  \ \{\tau^*<\infty\}.
\end{equation*}
If $\tau^*=\infty$ almost surely, 
then such a solution is called global.
\end{enumerate}

\end{definition}

\begin{remark}
This work focuses on classical solutions, i.e., solutions that are sufficiently regular for all terms to belong to $C^1(\K)$. Consequently, the parameter range for $s$ in the Sobolev space $H^s$ exceeds the optimal index $s = 3/2$ known for the deterministic case; see \cite{Guo-etal-2019-JDE, Tang-Shi-Liu-2015-MM}.
\end{remark}

\begin{remark}\label{Remark: t0 in R}
Due to the explicit dependence on $t$ in the damping coefficient $\lambda(t)$ and the noise coefficients $h_k(t,\cdot)$, it is necessary, when discussing asymptotic behaviour in this time-inhomogeneous system, to consider the Cauchy problem starting at an arbitrary fixed $t_0 \in \R$ (in this case the stochastic
integrals on $[t_0,t]$ are understood with respect to the increment processes
$W_k^{(t_0)}(t)=W_k(t)-W_k(t_0)$ and
$\widetilde W_k^{(t_0)}(t)=\widetilde W_k(t)-\widetilde W_k(t_0)$) rather than only at $t_0 = 0$
(see Remark~\ref{Remark:longtime_intro} below). Nevertheless, the particular choice of $t_0$ does not play an essential role in establishing the existence of solutions since only Brownian increments enter the
	stochastic integrals. For notational simplicity,  we shall therefore write $u_0 = u(t_0)$ or $u_0 = u(0)$ whenever the meaning of $t_0$ is clear from the context.
\end{remark}

Before delving into the concept of an evolution system of measures, we revisit some fundamental notions. Assuming the existence and uniqueness of the solution $u(t)$ to \eqref{Cauchy problem-Ito} with initial data $u(t_0) = u_0$ for any fixed time $t_0$ (as in Remark \ref{Remark: t0 in R}) and denoting $u(t,u_0)$ ($t_0 \le t$) as the state of $u$ at time $t$, the \textit{two-parameter} solution operator $\mathscr{S}_{t',t}$ of \eqref{Cauchy problem-Ito} is given by
\begin{align}
\mathscr{S}_{t',t}u(t',u_0) = u(t,u_0), \quad -\infty < t_0 \le t' \le t < \infty.\label{Solution operator}
\end{align}
Let $p$ be the transition kernel on $H^s$. The \textit{two-parameter} Markov propagator $\{\mathscr{P}_{t',t}\}_{t' \le t}$ on $f \in \mathscr{B}_B(H^s)$ is defined as follows:
\begin{align}\label{eq:markovsg}
\mathscr{P}_{t',t} f(\bullet) 
= \Ex f\big( \mathscr{S}_{t',t} u(t',\bullet)\big)
= \int_{H^s} f(v) p(t',t,\bullet,{\rm d}v).
\end{align}

This transition propagator is indeed Markov (due to the independent increments of the driving noise) and satisfies that for any $f \in \mathscr{B}_B(H^s)$ and $-\infty <  t_1 \le t_2 \le t_3 < \infty$,
\begin{align*}
\mathscr{P}_{t_1,t_2} \circ \mathscr{P}_{t_2,t_3} 
= \mathscr{P}_{t_1,t_3}.
\end{align*}
In fact, using the Markov property and the tower property of conditional expectation, we have  
\begin{align*}
	\mathscr{P}_{t_1,t_2} \circ \mathscr{P}_{t_2,t_3} f(u_0) 
= \Ex \Bigg[ \Ex \Big[ f \big( \mathscr{S}_{t_2,t_3} v \big) \Big] \Bigg|_{v = \mathscr{S}_{t_1,t_2}u(t_1,u_0)} \Bigg] = \Ex \Big[ \Ex \Big[ f \big( \mathscr{S}_{t_2,t_3} \circ \mathscr{S}_{t_1,t_2} u(t_1,u_0) \big) \Big| \mathcal{F}_{t_2} \Big] \Big] 
	= \mathscr{P}_{t_1,t_3} f(u_0).
\end{align*}
Let $\mathscr{P}^*_{t',t}$ be the adjoint map of $\mathscr{P}_{t',t}$. The above property of $\mathscr{P}_{t',t}$ implies 
$$\mathscr{P}^*_{t_2,t_3} \circ \mathscr{P}^*_{t_1,t_2} = \mathscr{P}^*_{t_1,t_3}, \quad  -\infty< t_1\leq t_2\leq t_3<\infty.$$

With this propagation property of $\mathscr{P}^*_{t',t}$, it is therefore reasonable to generalise the concept of the invariance property in the time-homogeneous case to the following sense (see \cite{DaPrato-Rockner-2008-Conference} for related results in a more general framework):
\begin{definition}\label{def:im_family}
A family of probability measures $\{\mu_t\}_{t \in \R}$ on the set $H^s$ is called \textit{an evolution system of measures} if
\begin{align}\label{eq:limitproperty}
\int_{H^{s}}[\mathscr{P}_{t',t}\phi](v)\mu_{t'}({\rm d}v)
=\int_{H^{s}}\phi(v)\mu_t({\rm d}v),
\quad \quad -\infty < t' \leq t < \infty,\quad \phi\in C_B(H^{s}).
\end{align}
\end{definition}

For more  details on evolution systems of measures, we refer to Remark \ref{Remark:longtime_intro} below.

Next, in order to state precisely our main results in Section \ref{sec:mainres_precise}, we will gather the hypotheses on the diffusion and damping coefficients as well as the hypotheses for the noises. We start by presenting the following straightforward assumptions on the diffusion and damping coefficients:
\begin{ManualHypo}
{\bf(${\mathbf H}_{\eps,\theta,\lambda}$)}\label{H-parameters}
For the diffusion term $\eps\Lambda^{2\theta}$ 
and  damping $\lambda(t)u$ in \eqref{Cauchy problem-Ito},  
we assume that $\eps\ge0$, $\theta\in(0,1]$ 
and $\lambda(\cdot)\in L^1_{\rm loc}(\R; [0,\infty))$.  
\end{ManualHypo}

Having detailed notation for discussing pseudo-differential operators in Section \ref{Section: Background}, we precisely state the assumptions we place upon the noise in \eqref{Cauchy problem-Ito}. We begin by specifying two classes of operators that are ``close" to skew-adjoint operators in the sense that an operator plus its adjoint is a zero-order operator. 

\begin{definition}\label{Ak class define}
Let  $\alpha\ge0$. We define
$$\mathbb{A}^{\alpha}\triangleq\Big\{\{a_k,\,\mathcal A_{k}\}_{k\ge1}\ :\ 
\text{the following conditions \ref{Ak:skew adjoint}\ 
and\  \ref{ak:coefficient} hold}\Big\},$$
where  \ref{ak:coefficient} and \ref{Ak:skew adjoint} are

\begin{enumerate}[label={$(\mathbb{A}_\arabic{enumi})$}]
\item\label{ak:coefficient} Either 
$$\{ a_k\}_{k\ge1}\subset H^\infty (\K;\R)\ \ \text{such that}\  \ 
\sum_{k=1}^\infty \norm{a_k }_{H^{r}} <\infty\ \ \text{for all}\ \ 
r\ge0,$$  or $$\{ a_k \}_{k\ge1}\in \ell^2;$$
\end{enumerate}
\begin{enumerate}[label={$(\mathbb{A}_2^\alpha)$}]

\item\label{Ak:skew adjoint} 
$\{\A_k\}_{k\ge1}\subset\OP\S^{\alpha}$, 
$\{\A_k+\A_k^*\}_{k\ge 1}\subset\OP\S^{0}$ are bounded.

\end{enumerate}

\end{definition}

\begin{definition}\label{Bk class define}
Let  $\beta\ge0$ and
$$\mathbb{B}^{\beta}\triangleq\Big\{\{b_k,\,\mathcal B_{k}\}_{k\ge1}\ :\ 
\text{the following conditions \ref{Bk:skew adjoint}\ 
and\  \ref{bk:coefficient} hold}\Big\},$$
where \ref{bk:coefficient} and \ref{Bk:skew adjoint} are
\begin{enumerate}[label={$(\mathbb{B}_\arabic{enumi})$}]
\item\label{bk:coefficient}  $\{ b_k \}_{k\ge1}\in \ell^2$.
\end{enumerate}

\begin{enumerate}[label={$(\mathbb{B}_2^\beta)$}]
\item\label{Bk:skew adjoint} 
$\{\B_k\}_{k\ge1}\subset\OP\S_0^{\beta}$;
$\{\B_k+\B_k^*\}_{k\ge 1}\subset\OP\S_0^{0}$ 
are bounded.
\end{enumerate}

\end{definition}

\begin{example}[Examples of classes 
\ref{Ak:skew adjoint} and \ref{Bk:skew adjoint}]

Let $k\ge 1$, $\alpha\in[0,1]$ and $\beta \ge0$. Assume that $\mathfrak{e}_k$ and $\mathfrak{h}_k$ are $\R$-valued  such that $\mathfrak{e}_k(\cdot)=-\mathfrak{e}_k(-\cdot)$ and $\mathfrak{h}_k(x,\cdot)=-\mathfrak{h}_k(x,-\cdot)$. Moreover, 
$$\{{\rm i}\mathfrak{h}_k\}_{k\ge1}\subset\S^0\ \ \text{and}\  \ \{{\rm i}\mathfrak{e}_k\}_{k\ge1}\subset \S_0^{0}\ \ \text{are bounded, respectively}.$$ 
Let $$\mathcal{H}_k\triangleq\OP({\rm i}\mathfrak{h}_k)\quad \mathcal{E}_k\triangleq \OP({\rm i}\mathfrak{e}_k).$$ 
It is easy to verify that  $\mathcal{E}_k$ is skew-adjoint and $\mathcal{H}_k$ satisfies $\mathcal{H}_k+\mathcal{H}_k^*\in \OP\S^{0}$. Then,
$\A_{k}\triangleq\mathcal{H}_k(\I-\Delta )^{\alpha/2}$ and  $\B_{k}\triangleq \mathcal{E}_k(\I-\Delta )^{\beta/2}$
satisfy \ref{Ak:skew adjoint} and  \ref{Bk:skew adjoint}, respectively.

\end{example}

Technically, even if $u\in H^\infty$, closing the $H^s$ 
estimate of $u$ is not obvious. Indeed, 
applying It\^o's formula to 
\eqref{Cauchy problem-Ito}$_1$, 
we encounter the following terms:
\begin{equation}\label{Qk-Cancel-example}
\sum _{k=1}^{\infty}\left \langle   \Q_k u,u\right \rangle _{H^{s}}
{\rm d}W_{k},\quad \sum _{k=1}^{\infty} 
\left \langle  \Q_k^{2} u,u\right \rangle _{H^{s}}{\rm d} t 
\quad \text{and} \quad
\sum _{k=1}^{\infty}\left
\langle  \Q_k u, \Q_k u\right \rangle _{H^{s}}{\rm d} t.
\end{equation}
If $\Q_k$ is of non-negative order, 
then all terms in \eqref{Qk-Cancel-example} 
are singular in the sense that they involve 
derivatives of order \textit{not less than} $s$. However, 
to close the $H^s$ estimate, they must 
be controlled using only $\|u\|^2_{H^s}$. 
Recently, such estimates for pseudo-differential 
operators $\Q_k=a_k\A_k$ or $\Q_k=b_k\B_k$ 
have been established in \cite{Tang-Wang-2022-arXiv,
Tang-2023-JFA} (see Lemmas 
\ref{Lemma-cancellation-1} and 
\ref{Lemma-cancellation-2-Xi} below).

In the following Hypothesis, 
we embed the necessary 
conditions to achieve this for $\Q_k$ containing
$a_k\A_k$ and $b_k\B_k$, namely:

\begin{ManualHypo}{\bf(${\mathbf H}_{\mathcal{Q}_k}$)}
\label{H-Q}
Assume that 
$$\Q_k = a_k\A_k + b_k\B_k\ \ \text{and}\ \ a_kb_k=0\ \ 
\text{for all}\ k\ge1,
$$ 
where $\{a_k,\,\mathcal A_{k}\}_{k\ge1}\in \mathbb{A}^\alpha$ 
with $\alpha\in[0,1]$ (cf. Definition \ref{Ak class define}), 
$\{b_k,\,\mathcal B_{k}\}_{k\ge1} \in\mathbb{B}^\beta$ 
with $\beta\ge0$ (cf. Definition \ref{Bk class define}).
\end{ManualHypo}

Next we state the hypotheses on 
$\sum_{k = 1}^\infty h_{k}(t,u)\d \widetilde{W}_{k}$. 
We require $\sum_{k = 1}^\infty h_k(t,\cdot)$ 
to be locally Lipschitz  in  $H^s$ and its growth to be 
bounded by a function containing a non-linear part 
in $\|\cdot\|_{\Wlip}$ and linear part in $\|\cdot\|_{H^s}$.

\begin{ManualHypo}{\bf(${\mathbf H}_{h_k}$)}\label{H-h}
Let $s \ge 1$.
$h_k\in \mathscr{B}(\mathbb R\times H^s; H^s)$, $k \in \N$. 
We assume there exist functions $K_1,K_2\in\mathscr{K}$ 
(defined in \eqref{SCRK})
such that for all $t\in\R$ and $u,v \in H^s$:
\begin{equation*}
\sum_{k=1}^\infty\|h_k(t,u)\|^2_{H^{s}}
\leq K_1(t,\|u\|_{\Wlip}) (1+\|u\|^2_{H^{s}}),\quad s\ge1,
\end{equation*}
\begin{equation*}
\sum_{k=1}^\infty \| h_k(t,u)- h_k(t, v)\|_{H^{s}}^2 
\leq K_2(t, \|u\|_{H^{s}}+\|v\|_{ H^{s}})\|u-v\|^2_{H^{s}},\quad s > 3/2.
\end{equation*}
\end{ManualHypo}

\begin{example}[Examples of $h_k$ satisfying Hypothesis \ref{H-h}]
Let $\mathcal{R}\in\OP\S_0^{-1}$ and $q(t)$ be a locally square integrable function in $t\in\R$. Then it is easy to see that
$$h_k(t,u)=2^{-k}q(t)\mathcal{R} \big[u^2+(\pd_{x} u)^2\big],\quad k\ge1$$
satisfies Hypothesis \ref{H-h}.

\end{example}

Under Hypotheses \ref{H-Q} and \ref{H-h}, 
we can show local-in-time existence of a unique 
solution $u$ to \eqref{Cauchy problem-Ito}. 
However, such a solution $u$ may blow up in finite time. 
Motivated by \cite{Tang-2023-JFA,Rohde-Tang-2021-NoDEA,
Alonso-Miao-Tang-2022-JDE}, we consider the 
regularisation effect from \textit{quickly growing} noise, 
which is described by the following Hypothesis \ref{H-h-large}:

\begin{ManualHypo}{\bf(${\mathbf H}^V_{h_k}$)}\label{H-h-large}

There are functions $V\in \mathscr{V}$ (cf. \eqref{SCRV}) and 
$g\in L_{\rm loc}^1(\mathbb{R};(0,\infty))$  
such that the following condition hold for $t\in\R,\ u\in H^s,\ s>3/2$:
\begin{align*} 
V'(\|u\|^2_{H^s})&\Big\{
(\Xi+2\Theta\|u\|_{\Wlip})\|u\|^2_{H^s}-2\lambda(t)\|u\|^2_{H^s}
+\sum_{k=1}^\infty\|h_k(t,u)\|^2_{H^s}\Big\} \\
&+2V''(\|u\|^2_{H^s})\left(\sum_{k=1}^\infty \bIP{h_k(t,u), u}_{H^s}^2 \right)
\leq g(t)V(\|u\|^2_{H^s}),
\end{align*}
where $\Xi=\Xi(\{a_k\}_{k\ge1},\{b_k\}_{k\ge1},s)$ 
and $\Theta=\Theta(s)$ are constants given 
in Lemmas \ref{Lemma-cancellation-2-Xi} 
and \ref{uux+F Hs inner product Te}, respectively.

\end{ManualHypo}

Indeed,  since $V''\leq0$, 
$2V''(\|u\|^2_{H^s})\big(\sum_{k=1}^\infty 
\bIP{h_k(t,u), u}_{H^s}^2 \big)$ 
can suppress the growth driven by other terms when the $\sum_{k = 1}^\infty h_k(t,\cdot)$ grows \textit{quickly enough}. In this way,  non-explosion can be ensured, and we refer to
\cite{Tang-Wang-2022-arXiv} for pertinent  
results in a more general  framework.

\begin{example}[Examples of $h_k$ satisfying  
Hypotheses \ref{H-h} and \ref{H-h-large}]
In order to substantiate the hypotheses  above, we give a non-trivial example of a noise $h$ 
and a Lyapunov function $V$ satisfying Hypotheses \ref{H-h} and \ref{H-h-large}.  

Let $\Psi(t,x)\in C\left(\mathbb{R}\times [0,\infty)\right)$ 
be locally bounded such that $\Psi(t,\cdot)$ is locally Lipschitz continuous. Besides, assume that $\Psi(t,x)\neq0$ for all $(t,x)\in \mathbb{R}\times [0,\infty)$, and for all $t\in\mathbb{R}$,
\begin{equation*}
\limsup_{x\rightarrow +\infty}\frac{2\Theta x}{\Psi^2(t,x)}<1,
\end{equation*}
where $\Theta$ is given in Lemma \ref{uux+F Hs inner product Te}.
Then, if $2\lambda(t)\ge 2\lambda_0\gg \Xi$ for all $t\in\mathbb{R}$, we have the following estimate for all $u\in H^s$ and $u\neq 0$:
\begin{align*}
&\frac{\left(\Xi-2\lambda(t) \right)\|u\|_{H^s}^2}{{\rm e}+\|u\|_{H^s}^2}
+\frac{2\Theta \|u\|_{\Wlip}\|u\|_{H^s}^2}{{\rm e}+\|u\|_{H^s}^2}\\
&+ \frac{\Psi^2(t,\|u\|_{\Wlip})\|u\|_{H^s}^2}{{\rm e}+\|u\|_{H^s}^2}-2\frac{\Psi^2(t,\|u\|_{\Wlip})\|u\|_{H^s}^4}{({\rm e}+\|u\|_{H^s}^2)^2} <0.
\end{align*}
This implies that
\begin{align*}
h_k(t,u) =\Psi(t,\|u\|_{\Wlip})\mathscr F^{-1}_{\xi\to x} \big(\mathbf{1}_{\{k-1<|\xi|\le k\}}  \mathscr F_{x\to\xi} u\big), \quad k\ge1
\end{align*}
satisfies
Hypotheses \ref{H-h} and \ref{H-h-large} with $V(x)=\log({\rm e}+x)$.

\end{example}

\subsection{Statements of  main results}
\label{sec:mainres_precise}
In this section, we record the precise 
statements of our main results, 
summarised in Section \ref{sec:mainres_summary}.
Before we do so, we recall the classes 
$\mathbb{A}^\alpha$ and $\mathbb{B}^\beta$ 
in Definitions \ref{Ak class define} and 
\ref{Bk class define}, respectively, in order to define 
a smoothness parameter:
\begin{equation}
\label{gamma0}
\gamma _{0}\triangleq
\begin{cases}
\max \Big\{ \alpha{\mathbf{1}}_{\big\{\sum _{k=1}^{\infty}\|a_{k}\|_{H^{s}}>0
\big\}},\
\beta{\mathbf{1}}_{\big\{\sum _{k=1}^{\infty}|b_{k}|^{2}>0\big\}} \Big\},
&\text{if}\ \{a_{k}\}_{k\ge 1}\subset H^{\infty}(\mathbb K;
\mathbb R),
\\[10pt]
\max \Big\{ \alpha{\mathbf{1}}_{\big\{\sum _{k=1}^{\infty}|a_{k}|^{2}>0
\big\}},\
\beta{\mathbf{1}}_{\big\{\sum _{k=1}^{\infty}|b_{k}|^{2}>0\big\}} \Big\},
&\text{if}\ \{a_{k}\}_{k\ge 1}\in l^{2}.
\end{cases}
\end{equation}

Our first result provides the local existence of a classical solution
and a suitable blow-up criterion.
\begin{theorem}\label{Thm:local-theory}
Let Hypotheses  \ref{H-Q}, \ref{H-h}  and \ref{H-parameters} hold. 
If $u_0$ is an $H^s$-valued $\mathcal{F}_0$-measurable 
random variable with 
$s>3/2+\max\{2\gamma_0,1,2\theta{\bf 1}_{\{\eps>0\}}\}$, 
where  $\gamma_0$ is given in \eqref{gamma0}, then
\begin{enumerate}[label={${\bf (\arabic*)}$}]
\item  \eqref{Cauchy problem-Ito} admits a unique 
maximal solution $(u,\tau^*)$ in the sense of 
Definitions \ref{pathwise solution definition};
\item $(u,\tau^*)$ defines a map 
$H^s\ni u_0\mapsto u(t)\in C([0,\tau^*);H^s)$ $\pas$, 
where $\tau^*$ does \textit{not} depend on $s$ 
and satisfies the following blow-up criterion
\begin{equation}
\label{blow-up criterion-1}
\displaystyle
{\bf 1}_{\left\{\limsup_{t\rightarrow \tau^*}\|u(t)\|_{H^s}=\infty\right\}}
={\bf 1}_{\left\{\limsup_{t\rightarrow \tau^*}\|u(t)\|_{\Wlip}
=\infty\right\}} \ {\rm on}\ \{\tau^*<\infty\}.
\end{equation}
\end{enumerate}

\end{theorem}

Exploiting the explicit blow-up criterion 
of Theorem \ref{Thm:local-theory}, it is 
possible to rule out blow-up with sufficient 
noise, (non-local) diffusion, and damping:
\begin{theorem}\label{Thm: global 1}
Let the conditions in 
Theorem \ref{Thm:local-theory} hold.
Assume that $K_1(\cdot,\cdot)\equiv c_0>0$ in Hypothesis  \ref{H-h}.
Then the following properties hold:
\begin{enumerate}[label={${\bf (\arabic*)}$}]
\item\label{T1:blow-up criterion}$($\textbf{Refined blow-up criterion}$)$ 
The blow-up criterion \eqref{blow-up criterion-1} can be improved to
\begin{equation}
\label{blow-up criterion-2}
\displaystyle
{\bf 1}_{\left\{\limsup_{t\rightarrow \tau^*}\|u(t)\|_{H^s}=\infty\right\}}
={\bf 1}_{\left\{\limsup_{t\rightarrow \tau^*}\int_{0}^{t}
\|\pd_x u(t')\|_{L^{\infty}}
\d t'=\infty\right\}} \ {\rm on}\ \{\tau^*<\infty\}.
\end{equation}

\item\label{T1:global} $($\textbf{Global existence}$)$  
If $\eps > 0$ and $\theta > \frac{1}{2}$, then $\mathbb{P}(\tau^* = \infty) = 1$.

\item\label{T1:asymptotics} $($\textbf{Asymptotics}$)$  
Let $\inf_{t\in\R}\lambda(t) \triangleq\lambda_0$. Recall that $\Xi$, defined in Lemma~\ref{Lemma-cancellation-2-Xi}, quantifies the strength of the noise term $\sum_{k = 1}^\infty \Q_k u(t) \circ{\rm d} W_k(t)$. 
For $\eps > 0$ and $\theta > \frac{1}{2}$, we have:
\begin{itemize}
	\item If $\lambda(t) \equiv \lambda_0 > \frac{\Xi + c_0}{2}$, then 
	\begin{align}
		\limsup_{t \to \infty} \Ex\left[\norm{u(t)}_{H^1}^2 \Big| \F_0 \right] 
		\leq \frac{c_0}{2\lambda_0 - \Xi - c_0} \quad \pas \label{limsup Ex u}
	\end{align}

\item If $\lambda_0 > \frac{\Xi + c_0}{2}$ and $\lim_{t\to\infty} \lambda(t) = \infty$, then
\begin{align}\label{lim Ex u=0}
\lim_{t \to \infty}\Ex\left[\norm{u(t)}_{H^1}^2\Big|\F_0\right]=0\quad \pas
\end{align} 
\end{itemize}

\end{enumerate}

\end{theorem}

Finally, exploiting large noise, we can 
establish the global regularity of the solutions in the absence 
of diffusion and exhibit an evolution system of measures:
\begin{theorem}\label{Thm: global 2}
In addition to the conditions in 
Theorem \ref{Thm:local-theory},  
if Hypothesis \ref{H-h-large} is also verified, 
then the following properties hold with $\sigma\in
\left(\frac{3}{2},{s-\max \{2\gamma_{0},1,2\theta{\bf 1}_{\{\eps>0\}}\}}\right)$:
\begin{enumerate}[label={${\bf (\arabic*)}$}]
\item\label{T2:global}$($\textbf{Global-in-time estimate}$)$ 
$\p(\tau^*=\infty)=1$, and 
\begin{equation}\label{Hs bound V}
\Ex\left[V(\norm{u(t)}_{H^s}^2)\Big|\F_0\right]  
\leq V(\| u_0\|^2_{H^s}) {\rm e}^{\int_0^t g(t')\d t'},
\quad t>0.
\end{equation}

\item\label{T2:stability}$($\textbf{Continuity of $u_0\mapsto u(t)$}$)$ 
If $\{u_{0,n}\}_{n\ge1}$ is a sequence of  
$\mathcal{F}_0$-measurable $H^s$-valued 
random variables such that 
$$
\sup_n\|u_{0,n}\|_{H^s}<\infty,\quad \lim_{n\to\infty}\norm{u_{0,n}-u_0}_{H^{\sigma}}=0\quad \pas,
$$ 
then the corresponding sequence of 
solutions $\{u_{n}\}_{n\ge1}$ satisfies
\begin{equation}\label{H-sigma stabilty}
\lim_{n\to\infty}
\Ex\left[1\land  \norm{u_{n}(t)-u(t)}^2_{H^{\sigma}}\Big|\F_0\right]=0, \quad t>0 \quad \pas
\end{equation}

\item\label{T2:invariant measure}
$($\textbf{An evolution system of measures}$)$ 
Moreover, if \eqref{Cauchy problem-Ito}  is considered in $\T$ and $g(\cdot)$ in Hypothesis \ref{H-h-large}  satisfies $g(\cdot) \in L^1(\R; (0,\infty))$, then there exists an evolution system of measures in the sense of Definition \ref{def:im_family}.

\end{enumerate}
\end{theorem}

\subsection{Remarks on the main difficulties and novelties}\label{Sect:remarks}

We assemble here several remarks on our results stated in the Section \ref{sec:mainres_precise}, highlighting key difficulties and novel approaches developed to overcome them. 

\begin{remark}
	The first point of emphasis is that, the noise structure $\sum_{k = 1}^\infty \Q_k u(t) \circ \mathrm{d} W_k(t) + \sum_{k = 1}^\infty h_k(t, u(t))  \mathrm{d} \widetilde{W}_k(t)$ in \eqref{Cauchy problem-Ito} encompasses both (Stratonovich) transport-type noise $\Q_k u \circ \mathrm{d} W_k$ (when $\Q_k$ reduces to a classical first-order differential operator) and (It\^o) multiplicative noise $h_k(t, u)  \mathrm{d} \widetilde{W}_k$.  This combined structure provides enhanced modelling flexibility through its incorporation of non-local pseudo-differential operators and singular interactions, facilitating the description of complicated non-local random interactions in physical systems. Furthermore, our analysis of such noise structures—to be elaborated in subsequent remarks—demonstrates significant potential for extension to many other models and generalization of prior results.

\end{remark}

\begin{remark}\label{Remark-singular}
Concerning Theorem~\ref{Thm:local-theory}, we note the following points:

\begin{itemize} 
 
 \item The first point of emphasis is the challenge posed by the singularity in \eqref{Cauchy problem-Ito}$_1$, specifically the inclusion of terms that lack invariance in $H^s$. If $u(t)$ is only known to be in $H^s$, directly applying It\^o's formula to $\|u(t)\|^2_{H^s}$ is problematic because the inner products $\langle u(t), u(t)\partial_{x}u(t) \rangle_{H^s}$, $\langle u(t), \Lambda^{2\theta} u(t) \rangle_{H^s}$, and $\langle u(t), \Q_k^2 u(t) \rangle_{H^s}$ are undefined. We overcome this issue by initially applying It\^o's formula to $\|J_n u(t)\|^2_{H^s}$, where $J_n$ is a mollifier operator defined in \eqref{Define Jn}, and then proceeding to the limit. This manoeuvre is employed repeatedly throughout the analysis, see, for example, \eqref{Ito to Jn u-1} and \eqref{Ito to Jn u-2}.

\item  Methodologically, our approach substantially
	extends previous methods for strong solutions in at least two key aspects:

\begin{itemize}
	\item We provide a unified framework for establishing local-in-time theory in both periodic ($x \in \mathbb{T}$) and whole-line ($x \in \mathbb{R}$) settings. For stochastic evolution systems on bounded domains like $\mathbb{T}$, the ``stochastic compactnes'' method—relying on tightness and Skorokhod's representation theorem—is commonly employed. See \cite{Alonso-Rohde-Tang-2021-JNLS, Tang-2018-SIMA, Holden-Karlsen-Pang-2023-DCDS, Galimberti-etal-2024-JDE} for applications to stochastic Camassa–Holm-type equations. Specifically, \cite{Alonso-Rohde-Tang-2021-JNLS} develops an abstract framework for transport-noise-driven systems, while \cite{Galimberti-etal-2024-JDE} provides recent refinements of the Skorokhod–Jakubowski theorem in quasi-Polish spaces. However, the non-compact embedding $H^s(\mathbb{R}) \subset H^\sigma(\mathbb{R})$ for $s > \sigma$ \textbf{obstructs} direct application of these methods on $\mathbb{R}$. Our work overcomes this issue through a direct approach, leveraging the recent framework for singular stochastic systems in general Hilbert spaces developed by the third author in \cite{Tang-Wang-2022-arXiv}.
	
	\item  In contrast to previous works on stochastic fluid models, our approach accommodates $H^s$-valued $\F_0$-measurable initial data $u_0$ \textbf{without} requiring moment conditions. This is achieved by employing conditional expectation $\mathbb{E}[\cdot|\F_0]$ in the proofs (see Lemma~\ref{Lemma-un convergence}). A significant advantage is that the blow-up criterion \eqref{blow-up criterion-1} emerges naturally from the existence proof, eliminating the need for separate arguments.  This substantially extends prior results where blow-up criteria were either absent or required a separate proof. See \cite{Crisan-Flandoli-Holm-2019-JNLS, GlattHoltz-Vicol-2014-AoP, Lang-Crisan-2022-SPDE, Chen-Duan-Gao-2023-AAP}.
\end{itemize}

\end{itemize}

\end{remark}

\begin{remark} Some remarks regarding global existence are in order:
	\begin{itemize}
		\item In Theorem~\ref{Thm: global 1}, when $\sum_{k=1}^\infty h_k(t,u)$ exhibits linear growth in $u$, we strengthen \eqref{blow-up criterion-1} to the following criterion:
		$$
		{\bf 1}_{\left\{\limsup_{t\rightarrow \tau^*}\|u(t)\|_{H^s}=\infty\right\}}
		={\bf 1}_{\left\{\limsup_{t\rightarrow \tau^*}\int_{0}^{t}
			\|\pd_x u(t')\|_{L^{\infty}}
			\d t'=\infty\right\}} \quad \text{on } \{\tau^*<\infty\}.
		$$ 
		Then, with suitable viscosity in the estimate of the $H^1$ norm, we can control this quantity, leading to global existence and demonstrating how damping affects the decay properties of $\|u(t)\|_{H^1}$. 
		\item For Theorem~\ref{Thm: global 2}, global regularity follows from direct estimates on $u(t)$ for all $t$ via the Lyapunov function (see \eqref{Hs bound V}). Notably, the strategy of employing rapidly growing (It\^o) multiplicative noise to control (potential) blow-up has been established in abstract frameworks for distribution-dependent settings \cite{Ren-Tang-Wang-2024-POTA} and singular stochastic evolution systems \cite{Tang-Wang-2022-arXiv}, with the latter covering stochastic incompressible Navier-Stokes and Euler equations. Building on the analysis of large noise in \cite{Tang-Wang-2022-arXiv,Ren-Tang-Wang-2024-POTA}, further examples and advances in long-time behaviour analysis are provided in \cite{Bagnara-Mario-Xu-2025-EJP,Tang-Wang-2024-CCM,Tang-Yang-2023-AIHP,Zhao-Li-Chen-2023-ZAMP,Chen-Miao-Shi-2023-JMP,Karlsen-Tang-Wang-2026-arXiv}.  In contrast to previous works, in \ref{T2:stability} of Theorem \ref{Thm: global 2},  rapidly growing noise additionally guarantees the continuity of the mapping $u_0 \mapsto u(t)$ for all $t$, a property essential for investigating the evolution system of measures—a topic which will be discussed in detail in the next remark.

	\end{itemize}

\end{remark}

\begin{remark}\label{Remark:longtime_intro} 

Regarding the evolution system of measures, we highlight several challenges related to the singularities and the time-inhomogeneous setting:
\begin{itemize} 

\item In the current time-inhomogeneous singular case, 
$\mathscr{P}_{t',t}$ deviates from the standard setup 
for the well-known Krylov--Bogoliubov argument. This deviation introduces some difficulties.

\begin{itemize} 
\item 
Since the diffusion can be \textit{degenerate} (i.e., $\varepsilon=0$), the first difficulty arises from the \textbf{mismatch} and \textbf{restriction} issues of the Markov propagator. As noted in Remark~\ref{Remark-singular}, the singular nature of the coefficients in \eqref{Cauchy problem-Ito}$_1$ poses an obstacle to establishing the Feller property for $\mathscr{P}_{t',t}: C_B(H^s) \to C_B(H^s)$. 
Recall the space  $X_{s_1,s_2}$ defined in \eqref{X-s1s2}.  Under the conditions of Theorem~\ref{Thm: global 2}, we only establish that for any $\varphi \in C_B(X_{s,\sigma})$ for some $\sigma<s$,  the restriction of $\mathscr{P}_{t',t}\varphi$ to closed $H^s$-ball is continuous with respect to the $H^\sigma$-topology (see \eqref{eq: restricted+mismatched feller} below). Note that 
the incompleteness of $X_{s,\sigma}$ implies that (weak) limits of measure sequences are \textbf{not a priori} guaranteed to remain in $\mathbf{P}(X_{s,\sigma})$ but may instead lie in  $\mathbf{P}(X_{\sigma,\sigma}) = \mathbf{P}(H^{\sigma})$. As far as we know, 
the mismatch  and restriction issues have not simultaneously been addressed in the literature, except  for the recent work \cite{Karlsen-Tang-Wang-2026-arXiv}. For instance, a  situation involving a discrepancy between the topology for the Feller property and the natural solution space was addressed in \cite[Section 3]{CotiZelati-GlattHoltz-Trivisa-2021-AMO} for the compressible 1D Navier--Stokes equations, in the context of a time-homogeneous Markov propagator.  However, their approach cannot handle the case when the Markov propagator has to be restricted. In this work, motivated by \cite{Karlsen-Tang-Wang-2026-arXiv}, we overcome both the \textit{mismatch} and \textit{restriction} issues of the Markov propagator. 
Specifically, we first identify a candidate for an evolution system of measures in $\mathbf{P}(H^\sigma)$. We then prove that this candidate actually lies in the smaller space $\mathbf{P}(X_{s,\sigma})$, with weak convergence upgradable to convergence in 
$\mathbf{P}(X_{s,\sigma})$.

\item The second difficulty arises from the absence of 
time homogeneity. In this case, by using the same method 
in proving Theorems \ref{Thm:local-theory} 
and \ref{Thm: global 2}, 
we obtain a global solution to  
\eqref{Cauchy problem-Ito} starting at any 
fixed $t_0 \in \mathbb{R}$ rather than $t_0 = 0$. 
The reason for this is that, in order to obtain \eqref{eq:limitproperty}, 
it is necessary to average backwards in time.
Assuming weak convergence, for an integer $n \ge 0$ and an initial measure $\mu$,
let $\nu_{-n}$ be a weak limit point as $T \to \infty$ of
\begin{equation}
\nu_{T,-n} =  \intbar_{-T}^{-n} \mathscr{P}^*_{t', -n} \mu \, \mathrm{d} t'.\label{approx-nu}
\end{equation}
Following \cite[Theorem 3.1]{DaPrato-Rockner-2008-Conference}, 
the evolution system of measures $\{\mu_t\}_{t \in \mathbb{R}}$ 
can be constructed by
\[
\mu_t \triangleq \mathscr{P}^*_{-n,t} \nu_{-n}, \qquad \nu_{-n} 
\triangleq \lim_{T \uparrow \infty} \nu_{T,-n},
\]
where the limit is understood as weak limit 
of the tight sequence $\{\nu_{T,-n}\}_{T > n}$ 
in the space $\mathbf{P}(X_{s,\sigma})$ (see \eqref{nu convergence X} below). 
As demonstrated in \eqref{eq: mu-t n-independence}, 
this definition of $\mu_t$ is independent of $n$. 
To appreciate the reason 
we average backwards to $t = -\infty$, 
we point the reader towards the property \eqref{eq:limitproperty}. 
If instead we averaged forward in time, 
defining the limit
\begin{equation*} 
\tilde{\nu}_n \triangleq 
\lim_{T \uparrow \infty}  \intbar_n^T \mathscr{P}^*_{n,t'} \mu \, \mathrm{d} t',
\end{equation*}
we see that this forward average integrates over the varying target time $t'$. This mixes states from different times, making it impossible to build $\tilde{\mu}_t$ satisfying the point-to-point transition required by \eqref{eq:limitproperty}.
\end{itemize}

\item In the time-homogeneous case, the operator $\mathscr{P}_{t',t}$, as defined in \eqref{eq:markovsg}, depends exclusively on the difference $t - t'$, thereby reducing to a one-parameter Markov semigroup $\{\mathscr{P}_{t}\}_{t \ge 0}$ with $t'$ set at $t'\equiv 0$. In the time-inhomogeneous case, the evolution system of measures $\{\mu_t\}_{t\in\R}$ broadens this aspect by emphasizing the transition from an \textit{arbitrary} time $t'$ to another \textit{arbitrary}  time $t$, rather than solely from the initial time $t_0$ to the present time $t$. This necessitates backward time averaging (in order to validate \eqref{eq:limitproperty}) and the adoption of the limit as $T\to\infty$ in \eqref{approx-nu} (refer to \eqref{eq:approx_measures} for the formal one), thus prompting a natural question:
\begin{center}
		\textit{``Which system's long-term dynamics does this measure system actually describe?"}
\end{center}
The answer is: it characterizes \textbf{not} a specific system but rather a hypothetical system with initial time $t_0\downarrow-\infty$. It is vital to emphasize that, despite the necessity of $t_0\downarrow-\infty$ for validating the limit of $T\to\infty$ (as in \eqref{eq:approx_measures}), there is \textbf{no} need to define the system's ``limit". This is justified by the important property of $n$-independence (see \eqref{eq: mu-t n-independence} below). Conceptually, one can imagine that the system starts at an arbitrary negative time, or simply $t_0=-\infty$, so that \eqref{eq:limitproperty} is satisfied for all $-\infty<t'\le t<\infty$.

\end{itemize}

\end{remark}

\section{Auxiliary results}\label{sec:aux_lemmas}

In this section, we collect some 
useful facts on mollification, commutators, operator bounds, and cancellation properties of the pseudo-differential noise. Recall that
$\OP$ is given in \eqref{OP define 1} and \eqref{OP define 2} on $\K=\R$ or $\K=\T$. For
$n\ge 1$, we define the mollifier on $\K$:
\begin{align}
J_{n} \triangleq{\mathrm{OP}}\big(j(\cdot /n)\big),\ \ n\ge 1,
\label{Define Jn}
\end{align}
where $j$ lies in the Schwartz space of
rapidly decreasing $C^{\infty}$ functions on $\mathbb R$ satisfying 
$0\leq j(y)\leq 1$ for all $y\in \mathbb R$ and $j(y)=1$ for any
$|y|\leq 1$. 

From the definitions of $\D^s$ and $J_n$, one can derive the following estimates:
\begin{lemma}\label{Lemma-Jn}
Let $\D^s$ and $J_n$ be defined 
in \eqref{Ds define} and \eqref{Define Jn}, respectively.
Then the following properties hold:
\begin{enumerate} 
\item[$(1)$] For all $\sigma\ge0$ and $n\ge1$ and $f\in H^\sigma$,   $\displaystyle\sup_n\|J_{n}\|_{\mathscr L(L^{\infty};L^{\infty})}<\infty$, $\displaystyle\sup_n\|J_{n}\|_{\mathscr L(H^{\sigma};H^{\sigma})}\leq 1$ and $\displaystyle\lim_{n\to\infty}  \|J_n f-f\|_{H^\sigma}=0$. Moreover, for all $r>s\ge0$,
\begin{align*}
\|{\mathbf I}-J_{n}\|_{\mathscr L(H^{r};H^{s})}\lesssim \,
\frac{1}{n^{r-s}},\quad 
\|J_{n}\|_{\mathscr L(H^{s};H^{r})}\sim \,  O(n^{r-s}).
\end{align*}
\item[$(2)$] For all $n\ge 1$ and $s\ge 0$,
$$
[\mathcal D^{s},J_{n}]=0,\ \ 
[\partial_x,J_{n}]=0,\ \
\left \langle J_{n}f, g\right \rangle _{L^{2}}=\left \langle f, J_{n}g
\right \rangle _{L^{2}}.
$$

\item[$(3)$] Let $f,g:\mathbb K\rightarrow \mathbb R$ satisfy
$g\in W^{1,\infty}$ and $f\in L^{2}$. Then, for some $C>0$,
\begin{align*}
\|[J_{n}, g]\partial_x f\|_{L^{2}} \leq C\|g\|_{W^{1,\infty }}\|f
\|_{L^{2}},\ \ n\ge 1.
\end{align*}
\end{enumerate} 

\end{lemma}

Now we recall some commutator estimates.

\begin{lemma}
[Proposition 3.6.A in \cite{Taylor-1991-book}]
\label{Taylor-commutator} 
Let $ s>0$ and $\P \in\OP\S^{s}$. 
Then for any $\sigma\ge0$, 
$g\in H^{s+\sigma}\cap W^{1,\infty},\ u\in H^{s-1+\sigma}\cap L^{\infty}$,
\begin{equation*}
\left\|\left[\P,g\I\right] u\right\|_{H^{\sigma}} 
\lesssim \|g\|_{W^{1,\infty}}\|u\|_{H^{s-1+\sigma}}
+ \|g\|_{H^{s+\sigma}}\|u\|_{L^{\infty}}.
\end{equation*}
\end{lemma}

\begin{lemma}[\cite{Tang-Zhao-Liu-2014-AA}]\label{F Lemma}
Recall the map $F(\cdot)$ defined in 
\eqref{Cauchy problem-Ito}.
Then we have
\begin{align*}
\|F(v)\|_{H^s}&\lesssim \|v\|_{W^{1,\infty}} \|v\|_{H^s},
\ \ v\in H^s,\ s>3/2,
\end{align*}
\begin{align*}
\|F(v_1)-F(v_2)\|_{H^s} 
&\lesssim\left(\|v_1\|_{H^s}+\|v_2\|_{H^s}\right)\|v_1-v_2\|_{H^s},
\ \   v_1,v_2\in H^s,\ s>3/2.
\end{align*}
\end{lemma}

\begin{lemma}\label{H1-conserve}
For all $u\in H^s$ with $s\ge2$ and $F(\cdot)$ 
defined in \eqref{Cauchy problem-Ito}, we have
\begin{align*}
\IP{u\,\pd_x u +F(u),u}_{H^1}=0.
\end{align*}
\end{lemma}

\begin{proof}
Direct computation shows  that for $u\in H^3$,
\begin{align*}
\IP{u\,\pd_x u +F(u),u}_{H^1}
=  \IP{\bk{\I -\pd^2_{x}}(u\,\pd_x u +F(u)),u}_{L^2} 
=  \IP{3u\,\pd_x u-2\pd_xu\,\pd^2_{x} u-u\,\pd^3_{x} u,u}_{L^2}=0.
\end{align*}
Let $2\le s<3$.
By density of $H^3\hookrightarrow H^{s}$, for $u\in H^s$, 
we can find $u_n\in H^3$ such that $u_n\to u$ in $H^s$, 
$F(u_n)\xrightarrow[]{n\to \infty} F(u)$, $u_n\,\pd_x u_n\xrightarrow[]{n\to \infty} u\,\pd_x u$ in $H^1$. 
Then we arrive at
\begin{align*}
0=\IP{u_n\,\pd_x u_n +F(u_n),u_n}_{H^1}\xrightarrow[]{n\to \infty} \IP{u\,\pd_x u +F(u),u}_{H^1},
\end{align*}
which is the desired result.
\end{proof}

\begin{lemma}
[\cite{Tang-Yang-2023-AIHP,
Miao-Rohde-Tang-2024-SPDE}]
\label{uux+F Hs inner product Te}
Let $s>3/2$. 
Let $F(\cdot)$ be given in \eqref{Cauchy problem-Ito}  
and $J_n$ be given in \eqref{Define Jn}. 
Then there is a constant $\Theta=\Theta(s)>0$ 
such that for all $n\ge1$,
\begin{equation*}
\left|\bIP{J_n  (uu_x), J_n  u}_{H^s}\right|
+\left|\bIP{J_n  F(u), J_n   u}_{H^s}\right|\leq \Theta \|u\|_{\Wlip}\|u\|^2_{H^s}.
\end{equation*}
\end{lemma}

\begin{lemma}\label{Fractional Laplacian}
Let $\Lambda^s$ be defined in \eqref{Lambda s define} and 
let $1/2<\sigma\leq 1$. Then the following inequality holds:
\begin{equation*} 
\norm{\partial_{x} f}_{L^{\infty}}
\lesssim \norm{\Lambda^{\sigma}f}_{H^{1}},\quad f\in H^{1+\sigma}.
\end{equation*}
\end{lemma}
\begin{proof}
By Sobolev embedding, 
for $0<\epsilon <\sigma-1/2$,  we have that
$$
\norm{\partial_{x} f}_{L^{\infty}}
\lesssim \norm{\partial_{x} f}_{H^{\frac{1}{2}+\epsilon}}
\lesssim \norm{\partial_{x} f}_{L^2}
+\norm{\Lambda^{3/2+\epsilon} f}_{L^2}. 
$$
Using the classical interpolation inequalities 
(cf. \cite[Proposition 1.32]{Bahouri-Chemin-Danchin-2011-book}), 
we have that
\begin{align}
\norm{\partial_{x} f}_{L^2}
\leq&\, \norm{\Lambda^\sigma f}_{L^2}^{\sigma}
\norm{\Lambda^{1+\sigma} f}_{L^2}^{1-\sigma}, \label{eq:interpol:1} \\
\big\|\Lambda^{3/2+\epsilon} f\big\|_{L^2}
\leq&\, \norm{\Lambda^\sigma f}_{L^2}^{\sigma-\frac{1}{2}-\epsilon} 
\norm{\Lambda^{1+\sigma} f}_{L^2}^{3/2+\epsilon-\sigma}, \label{eq:interpol:2}  
\end{align}
Finally, invoking Young's inequality 
$ab\leq \frac{1}{p}a^{p}+\frac{1}{q}b^{q}$ 
with $p=p_{1}=\frac{1}{\sigma},  q=q_{1}=\frac{1}{1-\sigma}$ 
in \eqref{eq:interpol:1} and 
$p=p_{2}=\frac{1}{\sigma-\frac{1}{2}-\epsilon},  
q=q_{2}=\frac{1}{\frac{3}{2}+\epsilon-\sigma}$ 
in \eqref{eq:interpol:2}, respectively, we find that 
\begin{align*}
&\norm{\pd_{x}f}_{L^2}
\leq \sigma \norm{\Lambda^\sigma f}_{L^2}
+ (1-\sigma) \norm{\Lambda^{1+\sigma} f}_{L^2},\\
\norm{\Lambda^{3/2+\epsilon} f}_{L^2}
&\leq (\sigma-1/2-\epsilon)
\norm{\Lambda^\sigma f}_{L^2}
+(3/2+\epsilon-\sigma)  
\norm{\Lambda^{1+\sigma} f}_{L^2}. 
\end{align*}
Taking into account that $0<\epsilon <\sigma-\frac{1}{2}\leq \frac{1}{2}$ 
and $\norm{\Lambda^{\sigma}f}^2_{H^{1}}
=\norm{\Lambda^{\sigma}f}^2_{L^2}
+\norm{\Lambda^{1+\sigma}f}^2_{L^2}$, 
the result follows.
\end{proof}

Finally, we state the following results on compositions of mollifications and pseudo-differential operators, showing that the terms in \eqref{Qk-Cancel-example} interact to yield cancellations uniformly in the mollifier.

\begin{lemma}[\cite{Tang-Wang-2022-arXiv,
Tang-2023-JFA}]\label{Lemma-cancellation-1}
Assume that $(a_k,\A_k)\in \mathbb{A}^\alpha$ with $\alpha\in[0,1]$ 
$($see Definition \ref{Ak class define}$)$ and 
$(b_k,\B_k)\in\mathbb{B}^\beta$ with $\beta\ge0$ 
$($cf. Definition \ref{Bk class define}$)$. 
Let $J_n$ be defined in \eqref{Define Jn} and 
$\gamma_0$ be given in \eqref{gamma0}.
Then for $\eta\ge0$, there is a constant 
$$
C=C(\{a_k\}_{k\ge1},\{b_k\}_{k\ge1},\eta)>0
$$
such that
\begin{equation*}
\sup_{n\ge 1} \sum_{k=1}^\infty 
\left( \IP{J_n (a_k\A_k)J_n f,\, f}_{H^\eta}^2 
+ \IP{J_n (b_k\B_k)J_n f,\, f}_{H^\eta}^2 \right)
\le C\|f\|_{H^{\eta}}^4,\ \ f\in  H^{\eta},
\end{equation*}
and
\begin{equation*}
\sup_{n\ge 1} \sum_{k=1}^\infty 
\left(\IP{J_n (a_k\A_k) f,\, J_n f}_{H^\eta}^2 
+\IP{J_n (b_k\B_k) f,\, J_n f}_{H^\eta}^2 \right)
\le C\|f\|_{H^{\eta}}^4,\ \ f\in  H^{\eta\lor \gamma_{0}}.
\end{equation*}
\end{lemma}

\begin{lemma}
[\cite{Tang-Wang-2022-arXiv,
Tang-2023-JFA}]\label{Lemma-cancellation-2-Xi}
Assume that $(a_k,\A_k)\in \mathbb{A}^\alpha$ with $\alpha\in[0,1]$ 
$($see Definition \ref{Ak class define}$)$ and 
$(b_k,\B_k)\in\mathbb{B}^\beta$ with $\beta\ge0$ 
$($cf. Definition \ref{Bk class define}$)$. 
Let $J_n$ be defined in \eqref{Define Jn} 
and $\gamma_0$ be given in \eqref{gamma0}.
Then for $\eta\ge1$, there is a constant
$$
\Xi={\Xi}(\{a_k\}_{k\ge1},\{b_k\}_{k\ge1},\eta)>0
$$
such that for $ f\in H^{\eta}$,
\begin{align*}
\sup_{n\ge 1}\, &\sum_{k=1}^{\infty}
\Big|\IP{J_n^3(a_k\A_k)^2 J_n f,f}_{H^\eta}
+\IP{J_n (a_k\A_k) J_n f,J_n (a_k\A_k) J_n f}_{H^\eta}\Big|\\
&+ \sup_{n\ge 1}\sum_{k=1}^{\infty}
\Big|\IP{J_n^3 (b_k\B_k)^2J_n f,f}_{H^\eta}
+\IP{J_n (b_k\B_k)J_n f,J_n (b_k\B_k)J_n f}_{H^\eta}\Big| 
\leq \,
\Xi \norm{f}^2_{H^\eta},
\end{align*}
and for $f\in H^{\eta\lor 2\gamma_{0}}$,
\begin{align*}
\sup_{n\ge 1}\sum_{k=1}^{\infty}\
&\Big|\IP{J_n(a_k\A_k)^2 f,J_nf}_{H^\eta}
+\IP{J_na_k\A_k f,J_n a_k\A_k f}_{H^\eta}\Big|\\
&+ \sup_{n\ge 1}\sum_{k=1}^{\infty}
\Big|\IP{J_n(b_k\B_k)^2 f,J_n f}_{H^\eta}
+\IP{J_nb_k\B_k f,J_nb_k\B_k f}_{H^\eta}\Big| 
\leq \,
\Xi \norm{f}^2_{H^\eta}.
\end{align*}
\end{lemma}

To conclude this section, we recall some facts regarding the weak convergence of probability measures. Let $(X,\rho)$ be an arbitrary metric space.  A family of probability measures $\mathcal{M}\subset \mathbf{P}(X)$ is tight if for each $\epsilon>0$, there is  a compact set $K_\epsilon\subset X$ such that 
$$\inf_{\mu\in\mathcal{M}}\mu(K_\epsilon)>1-\epsilon.$$

A sequence $\{\mu_n\}_{n\ge1}\subset \mathbf{P}(X)$ is called to converge weakly to $\mu\in\mathbf{P}(X)$ if 
$$\lim_{n\to\infty}\int_Xf(v)\mu_{n}({\rm d}v)=\int_Xf(v)\mu({\rm d}v),\quad f\in C_B(X).$$

\begin{lemma}[\cite{Ethier-Kurtz-1986-Book,Billingsley-1999-Book}] \label{convergence of measures lemma}
Let $(E,d_E)$ be a separable metric space and let $(\mathbf{P}(E),\widetilde d_E)$  be the corresponding space of probability measures equipped with the Prokhorov metric $\widetilde d_E$. If $\mathscr{Q}\subset \mathbf{P}(E)$ is tight, then it is relatively sequentially compact; that is, every sequence in $\mathscr{Q}$ contains a weakly convergent subsequence with limit in $ \mathbf{P}(E)$.
Moreover, for any sequence $\{\mu_n\}_{n\ge1}\subset \mathbf{P}(E)$ and $\mu\in \mathbf{P}(E)$, the following conditions are equivalent: 
\begin{enumerate}[leftmargin=0.79cm,label={${\bf (\alph*)}$}]\setlength\itemsep{0.2em}
	\item\label{measure conver 2.1} $\widetilde d_E(\mu_n,\mu)\to0$ as $n\to\infty$;
	\item\label{measure conver 2.2} $\mu_n$ converges weakly to $\mu$ as $n\to\infty$;
	\item For all $f\in C_{B,U}(E),\, \lim_{n\to\infty}\int_E f\,d\mu_n = \int_E f\,d\mu$;
	\item\label{measure conver 2.4} For all open subset $G\subset E,\, \mu(G) \le \liminf_{n\to\infty}\mu_n(G)$.
\end{enumerate}

\end{lemma}

\section{Local pathwise classical solution}\label{sec:localtheory}

We now present the proof for Theorem \ref{Thm:local-theory}.  Primarily, the proof draws upon the recent work by the third author in \cite{Tang-Wang-2022-arXiv}, which provides a novel framework for handling singular stochastic evolution systems in general Hilbert spaces. Despite this, given that some parameters in the proof will continue to be employed in subsequent sections, we opt to write down the comprehensive proof to ensure both completeness and consistency, as opposed to merely verifying the conditions from \cite{Tang-Wang-2022-arXiv}.

\subsection{Step 1: the  approximation scheme}
Given that either $a_{k}=0$ or
$b_{k}=0$ for all $k\ge 1$, we can express
\begin{equation*}
\sum _{k=1}^{\infty}\Q_{k} u(t) \d W_{k}(t)=
\sum _{k=1}^{\infty}a_k\A_{k} u(t) \d\overline W_{k}(t)+
\sum _{k=1}^{\infty} b_k\B_{k} u(t) \d\widehat W_{k}(t)
\end{equation*}
where $\big\{\overline W_{k}(t),  \widehat W_{k}(t)\big\}_{k\ge 1}$ denote a family of independent standard 1-dimensional Brownian motions, which are
also independent of $\{\widetilde{W}_{k}(t)\}_{k\ge 1}$.
Let $\U$  be a separable Hilbert space with 
a complete orthonormal basis $\{e_k\}_{k\ge 1}$. 
Define a cylindrical 
Brownian motion $\W$  on $\U$ as
\begin{equation*}
\W(t)\triangleq\sum_{k=1}^\infty \left(\ol W_k(t) e_{3k-2} +\hh W_k(t) e_{3k-1}+\widetilde{W}_k(t) e_{3k}\right)
\end{equation*} 
and define
the quantity $H(t,u)$ unifying the noise coefficient as follows:
\begin{equation*} 
H(t,u)e_{3k-2} \triangleq a_k\A_ku,\ \  
H(t,u)e_{3k-1} \triangleq  b_k\B_k  u,\ \ 
H(t,u)e_{3k}\triangleq h_k(t,u),\ \ k\ge 1.
\end{equation*} 
Then we can infer from \eqref{OP boud Hs},
Hypotheses \ref{H-Q} and \ref{H-h} that
$
H:\R\times H^s\ni(t,u)\mapsto H(t,u)\in\LL_{2}(\U;H^{s-\gamma_0})
$
is measurable.  Now {\eqref{Cauchy problem-Ito}$_1$} can be rewritten as
\begin{align*}
\d u(t)+\Big[\eps\Lambda^{2\theta} u +\lambda(t)\, u+u\,\pd_x u +F(u)\Big](t) \d t 
= \frac{1}{2}\sum _{k=1}^{\infty}\Big[ (a_k\A_{k})^{2} u +(b_k\B_{k})^{2} u\Big](t)
\d t+H(t,u(t))\d \W(t).
\end{align*}

As noted in Remark \ref{Remark-singular}, the aforementioned equation is singular in $H^s$ due to the non-invariance of its drift and noise coefficients in $H^s$. Additionally, the equation is non-linear, and in general, the expectation of products cannot be separated. Consequently, we proceed to mollify the equation and truncate the non-linearity.

For any $R\ge1$, we take a cut-off function 
$\chi_R\in C^{\infty}([0,\infty);[0,1])$ such that 
$\chi_R(r)=1\ \text{for}\ |r|\le R,\ 
\text{and}\ \chi_R(r)=0\ \text{for}\ r>2R.$ 
Let $\{J_n\}_{n\ge 1}$ be the mollifier
given in \eqref{Define Jn}.  Define 
\begin{equation*}
H_{n,R}(t,v):\mathbb{R}\times H^s\ni(t,v)
\mapsto H_{n,R}(t,v)\in\LL_{2}(\U;H^{s}),\quad n\ge1,\ R\ge1,\ s>3/2
\end{equation*} 
such that for all $k\ge1$ and $v\in H^s$,
\begin{equation*}
H_{n,R}(t,v)e_{3k-2} \triangleq J_n(a_k\A_k) J_n v,\ \
 H_{n,R}(t,v)e_{3k-1} \triangleq J_n(b_k\B_k) J_nv,\ \
 H_{n,R}(t,v)e_{3k}\triangleq \chi_R\big(\|v-u_0\|_{\Wlip}\big)h_k(t,v).
\end{equation*} 
For the other terms of the equation, we set
\begin{equation*}
G_{n,R}(v)\triangleq\chi^2_R\big(\|v-u_0\|_{\Wlip}\big)
\Big[J_n[J_n v\,\pd_xJ_n v] +F(v)\Big],\ \ 
v\in H^s,\ s>3/2,
\end{equation*} 
where $F(u)$ is the non-linear, non-local part 
of the equation given in \eqref{F define}. For $n\ge1,\ R\ge1$, 
we consider the following approximation 
of \eqref{Cauchy problem-Ito}.
\begin{equation}\label{approximation cut eq}
\left\{\begin{aligned}
\d u(t)&+\Big[\eps J_n\Lambda^{2\theta}J_n u(t) 
+\lambda(t)\, u(t)+G_{n,R}(u(t))\Big]\d t\\
&= \frac{1}{2}\sum _{k=1}^{\infty}
\left [ J_n^3(a_k\A_{k})^{2} J_n u(t) 
+J_n^3 (b_k\B_{k})^{2} J_n u(t)\right ]
\d t+H_{n,R}(t,u(t))\d \W(t),\\
u(0)&=u_0.
\end{aligned}\right.
\end{equation}

Let $s>3/2+\max\{2\gamma_0,1,2\theta{\bf 1}_{\{\eps>0\}}\}$. 
Now we state some estimates on the coefficients in 
\eqref{approximation cut eq}$_1$, which will be used later. 
By a direct calculation, we have
\begin{equation}
2\IP{\eps J_n^2\Lambda^{2\theta} v +\lambda(t)\, v,v}_{H^s}
=2\eps\|\Lambda^{\theta} J_n v\|^2_{H^s}
+2\lambda(t)\|v\|^2_{H^s} \ge0,\ \ v\in H^s.\label{Hs estimate-1}
\end{equation}
And by Lemmas \ref{Taylor-commutator} and \ref{F Lemma}, 
\begin{equation}
2\IP{G_{n,R}(v), v}_{H^s}
\lesssim (\|u_0\|_{\Wlip}+2R)\|v\|^2_{H^s},\ \ v\in H^s.
\label{Hs estimate-2}
\end{equation}
By Hypotheses \ref{H-Q} and \ref{H-h}, 
Lemmas \ref{Lemma-cancellation-1} and \ref{Lemma-cancellation-2-Xi}, 
we find a function $K\in\mathscr K$ such that for all $v\in H^s$, 
$n\ge1$ and $t\ge 0$,
\begin{align}
&\sum_{k=1}^{\infty}\bIP{\left[ J_n^3(a_k\A_{k})^{2} J_n v 
+J_n^3 (b_k\B_{k})^{2} J_n v\right ], v}_{H^s} +\|H_{n,R}(t,v)\|_{\LL_2(\U;H^s)}^2 \nonumber\\
= \ & \chi^2_R\big(\|v-u_0\|_{\Wlip}\big)
\left(\sum_{k=1}^{\infty}\|h_k(t,v)\|_{H^{s}}^2\right)\notag\\
&+\sum_{k=1}^{\infty}\IP{J_n^3(a_k\A_k)^2J_nv,v}_{H^{s}}
+ \sum_{k=1}^{\infty} \|J_n (a_k\A_k) J_nv\|_{H^{s}}^2\nonumber\\
&+\sum_{k=1}^{\infty}\IP{J_n^3(b_k\B_k)^2J_nv,v}_{H^{s}}
+\sum_{k=1}^{\infty} \|J_n (b_k\B_k) J_nv\|_{H^{s}}^2\nonumber\\
\le\ & K(t,\|u_0\|_{\Wlip}+2R) (1+\|v\|_{H^s}^2),\label{Hs estimate-3}
\end{align} 
and 
\begin{align} 
&\sum_{k=1}^\infty\IP{H_{n,R}(t,v) e_k, v}_{H^s}^2 \notag\\
=\ &\sum_{k=1}^\infty\IP{J_n (a_k\A_k) J_nv,v}_{H^{s}}^2 
+\sum_{k=1}^\infty\IP{J_n (b_k\B_k) J_nv,v}_{H^{s}}^2 +\chi^2_R\big(\|v-u_0\|_{\Wlip}\big)
\sum_{k=1}^\infty\IP{h_k(t,v),v}_{H^{s}}^2\notag\\
\le \ &  K(t,\|u_0\|_{\Wlip}+2R) \big(1+ \|v\|_{H^s}^4\big).
\label{Hs estimate-4}
\end{align}

These estimates provide us with global classical solutions
to the approximate system \eqref{approximation cut eq}:
\begin{lemma}\label{Lemma-un T estimates} Under the conditions 
in Theorem \ref{Thm:local-theory}, for any $R,n\ge1$, 
\eqref{approximation cut eq} has a unique global solution 
$u_n^{(R)}\in C([0,\infty);H^s)$ $\pas$ Moreover, 
there is a function $Q: [0,\infty)\times [0,\infty)\to (0,\infty)$ 
increasing in both variables such that for all $R\ge1,\ T>0$,
\begin{align}\label{eq:unif_approx_gwp}
\sup_{n\ge1}\Ex
\left[\sup_{t\in[0,T]}\big\|u_n^{(R)}(t)\big\|^2_{H^s}\bigg|\F_0\right] 
\leq Q(T,2R+\|u_0\|_{\Wlip})(1+\|u_0\|_{H^s}^2).
\end{align} 
\end{lemma}

\begin{proof}
For any $n,R\ge1$, all the coefficients of 
\eqref{approximation cut eq} 
are locally Lipschitz continuous in $u\in H^s$ locally uniform in $t$. 
Consequently, for any deterministic initial data, 
\eqref{approximation cut eq} 
has a unique solution, which is continuous in $H^s$  \cite{Wang-2013-Book}.  
Since $\F_0$ is independent of the equation, 
\eqref{approximation cut eq} also admits a unique 
solution $u_n^{(R)}(t)$  with $u_n^{(R)}(0)=u_0$,  
and $u_n^{(R)}(t)\in C([0,\tau_n(R));H^s)$ with 
$$\tau_n(R)
\triangleq\lim_{N\to\infty}\tau_{n,N}(R),\ \ 
\tau_{n,N}(R)
\triangleq\inf\Big\{t\ge 0: \big\|u_n^{(R)}(t)\big\|_{H^s}\ge N\Big\},\quad N\ge 1.$$

Therefore, for any $T>0$, by using It\^o's formula with  
$v=u_n^{(R)}(t)$ in \eqref{Hs estimate-1}, \eqref{Hs estimate-2}, 
\eqref{Hs estimate-3}, and \eqref{Hs estimate-4}, 
and then employing the BDG inequality, 
we find constants $c_1,c_2>0$  such that 
for any $t\in [0,T]$ and $N\ge 1$, 
\begin{align*}
&\Ex\left[
\sup_{t'\in[0,t\land\tau_{n,N}(R)]}\big\|u_n^{(R)}(t')\big\|^2_{H^s}\bigg|\F_0\right]- \|u_0\|^2_{H^s}\\
\le \ &
c_1 \Ex\bigg[\bigg(\int_0^{t\land\tau_{n,N}(R)} K(t',\|u_0\|_{\Wlip}+2R) 
\Big(1+\big\|u_n^{(R)}(t')\big\|_{H^s}^4\Big)\d t'\bigg)^{\frac 1 2}\bigg|\F_0\bigg]\\  
 &+ c_1 \Ex\bigg[\int_0^{t\land\tau_{n,N}(R)} K(t',\|u_0\|_{\Wlip}+2R) \left(1+\big\|u_n^{(R)}(t')\big\|_{H^s}^2\right)\d t'\bigg|\F_0\bigg]\\
\le \ &  \frac12 \Ex\Big[\sup_{t'\in [0,t\land\tau_{n,N}(R)]} \big\|u_n^{(R)}(t')\big\|_{H^s}^2\Big|\F_0\Big] 
+ c_2\int_0^{t} K(t',\|u_0\|_{\Wlip}+2R)\d t'\\
 &+ c_2 \int_0^t K(t',\|u_0\|_{\Wlip}+2R) \Ex\left[\sup_{r\in [0,t'\land \tau_{n,N}(R)]} \|u_n^{(R)}(r)\|_{H^s}^2\Big|\F_0\right]\d t'.
\end{align*}
By Gr\"onwall's inequality, there exists a 
function $Q: [0,\infty)\times [0,\infty)\to (0,\infty)$, 
increasing in both variables, such that for all $n,N\ge 1$,
\begin{align} 
\Ex\left[
\sup_{t\in[0,T\land\tau_{n,N}(R)]}\big\|u_n^{(R)}(t)\big\|^2_{H^s}\bigg|\F_0\right] 
\le\, Q(T,2R+\|u_0\|_{\Wlip})(1+\|u_0\|_{H^s}^2).\label{EXN}
\end{align} 
This implies that for all $n,N\ge1$,
\begin{align*}
\p(\tau_{n,N}(R)<T|\F_0)
\le  \frac{Q(T,2R+\|u_0\|_{\Wlip})(1+\|u_0\|_{H^s}^2)}{N^2},
\end{align*} 
so that $\tau_n(R)\triangleq\lim_{N\to\infty}\tau_{n,N}(R)$ satisfies
$
\p(\tau_n(R)<T|\F_0)\le \lim_{N\to\infty} \p(\tau_{n,N}(R)<T|\F_0) =0.
$
Hence, $\p(\tau_n(R)\ge T)= \Ex[\p(\tau_n(R)\ge T|\F_0)]=1$ 
for all $T>0$, which means  $\p(\tau_n(R)=\infty)=1$, 
i.e., global existence. 
By letting $N\to\infty$ in \eqref{EXN}, we obtain the 
estimate \eqref{eq:unif_approx_gwp}.
\end{proof}

\subsection{Step 2: convergence of the approximation solutions}

\begin{lemma}\label{Lemma-un convergence} 
Let $R\ge1$ and $\{u_n^{(R)}\}_{n\ge1}$ be the sequence 
of approximate solutions given in Lemma \ref{Lemma-un T estimates}. Let  $\sigma\in\left(\frac{3}{2},{s-\max \{2\gamma_{0},1,2\theta{\bf 1}_{\{\eps>0\}}\}}\right)$.  
Under the conditions in Theorem \ref{Thm:local-theory}, 
the following properties hold:
\begin{enumerate}[label={\bf (\arabic*)},leftmargin=*]
\item \label{Global cut-off convergence} Define
$$\tau_N^{n,l}(R)\triangleq N\land \inf\big\{t\ge 0: \big\|u_n^{(R)}(t)\big\|_{H^s}\lor \big\|u_l^{(R)}(t)\big\|_{H^s}\ge N\big\}
,\quad n,l\ge 1, \ N\ge1.$$
Then $\pas$,
\begin{equation} \label{Cauchy in E-M}
\lim_{n\rightarrow\infty}\sup_{l\ge n}
\Ex\bigg[\sup_{t\in[0,\tau_N^{n,l}(R)]} \big\|u_n^{(R)}(t)-u_l^{(R)}(t)\big\|^2_{H^\sigma}\bigg|\F_0\bigg]
=0,\quad N>0.
\end{equation}

\item \label{Convergence of u-n} There is an 
$\mathcal{F}_t$-progressively measurable 
$H^s$-valued process 
$\big(u^{(R)}(t)\big)_{t\ge0}$  and a subsequence of $\{u_n^{(R)}\}$  
$($still labeled as $\{u_n^{(R)}\}$ for simplicity$)$ 
such that $\pas$,
\begin{equation}\label{Xn to X}
u_n^{(R)}\xrightarrow[]{n\to \infty}u^{(R)} \ {\rm in}\ 
C([0,T];H^\sigma),\quad T>0.
\end{equation}
Furthermore, for $Q(\cdot,\cdot)$ given in 
Lemma \ref{Lemma-un T estimates}, we have that $\pas$,
\begin{equation}\label{X L2 bound}
\Ex\Big[\sup_{t\in [0,T]} \|u^{(R)}(t)\|_{H^s}^2\Big|\F_0\Big]\le Q(T,2R+\|u_0\|_{H^s})(1+\|u_0\|^2_{H^s}),\ \ T>0.
\end{equation} 
\end{enumerate} 
\end{lemma}
\begin{proof}
\ref{Global cut-off convergence} 
Set $v_{n,l}^{(R)} \triangleq u_n^{(R)}-u_l^{(R)}$ for $n,l\ge 1$. 
By  It\^{o}'s formula,
\begin{align*}
{\rm d}\big\|v_{n,l}^{(R)}(t)\big\|_{H^\sigma}^2
= \Big(\sum_{i=1}^3 A_{i}(t)
+\sum_{k=1}^\infty \sum_{i=4}^8A_{i,k}(t)\Big)\d t +2\IP{
\Big(H_{n,R}(t,u_n^{(R)})-H_{l,R}(t,u_l^{(R)})\Big) {\rm d}\W(t),
v_{n,l}^{(R)}
}_{H^\sigma},
\end{align*}
where 
\begin{align*} 
A_{1}(t) &\triangleq -2\eps \IP{J_n \Lambda^{2\theta}J_n u_n^{(R)}(t)-J_l \Lambda^{2\theta}J_l u_l^{(R)}(t),\ v_{n,l}^{(R)}(t) }_{H^{\sigma}},\\
A_{2}(t) &\triangleq -2\lambda(t) \IP{v_{n,l}^{(R)}(t),\ v_{n,l}^{(R)}(t) }_{H^{\sigma}},\\
A_{3}(t)&\triangleq -2\IP{G_{n,R}\big(u_n^{(R)}(t)\big)-G_{l,R}\big(u_l^{(R)}(t)\big),\ v_{n,l}^{(R)}(t) }_{H^{\sigma}},\\
A_{4,k}(t)&\triangleq  \IP{J_n^3(a_k\A_k)^2J_n u_n^{(R)}(t)-J_l^3(a_k\A_k)^2J_l u_l^{(R)}(t),\ v_{n,l}^{(R)}(t) }_{H^{\sigma}},\\ 
A_{5,k}(t)&\triangleq \IP{J_n^3(b_k\B_k)^2J_n u_n^{(R)}(t)-J_l^3(b_k\B_k)^2J_l u_l^{(R)}(t),\ v_{n,l}^{(R)}(t) }_{H^{\sigma}},\\ 
A_{6,k}(t)& \triangleq \left\| J_n(a_k\A_k)J_n u_n^{(R)}(t)-J_l(a_k\A_k)J_l u_l^{(R)}(t)\right\|^2_{H^{\sigma}},\\  
A_{7,k}(t)&\triangleq \left\| J_n(b_k\B_k)J_n u_n^{(R)}(t)-J_l(b_k\B_k)J_l u_l^{(R)}(t)\right\|^2_{H^{\sigma}},
\end{align*} 
and
\begin{align*} 
A_{8,k}(t) 
\triangleq  \Big\| \chi_R\big(\|u_n^{(R)}(t)-u_0\|_{\Wlip}\big) h_k\big(t,u_n^{(R)}(t)\big)- \chi_R\big(\|u_l^{(R)}(t)-u_0\|_{\Wlip}\big) h_k\big(t,u_l^{(R)}(t)\big)\Big\|^2_{H^{\sigma}}.  
\end{align*} 
\textbf{Claim:} There exists a function $\hh K\in\mathscr{K}$ 
(cf. \eqref{SCRK}) and a function 
$\varrho: \mathbb N\times\mathbb N \rightarrow[0,\infty)$ 
with $\displaystyle \lim_{n,l\to\infty} 
\varrho_{n,l}=0$ such that
\begin{align}
 \sum_{i=1}^3 A_{i}(t)
+\sum_{k=1}^\infty \sum_{i=4}^8A_{i,k}(t) 
\leq     \hh K\left(t,\big\|u_n^{(R)}(t)\big\|_{H^s}
+\big\|u_l^{(R)}(t)\big\|_{H^s}\right)(\varrho_{n,l}
+\big\|v_{n,l}^{(R)}(t)\big\|^2_{H^\sigma}),
\label{Check convergence-1} 
\end{align}
and
\begin{align}
\sum_{k=1}^\infty& 
\IP{\{H_{n,R}(t,u_n^{(R)}(t))-H_{l,R}(t,u_l^{(R)}(t))\}e_k, v_{n,l}^{(R)}(t)}_{H^\sigma}^2 \nonumber\\                               
&\le  \hh K\left(t,\big\|u_n^{(R)}(t)\big\|_{H^s}
+\big\|u_l^{(R)}(t)\big\|_{H^s}\right)\big\|v_{n,l}^{(R)}(t)\big\|^2_{H^\sigma}
\left(\varrho_{n,l}+\big\|v_{n,l}^{(R)}(t)\big\|^2_{H^\sigma}\right),
\label{Check convergence-2} 
\end{align}

The terms $A_{i,k}$ $(i=4,5)$ arise from the 
Stratonovich-to-It\^o conversion, and the 
terms $A_{j,k}$ $(j=6,7,8)$ arise from the It\^o correction. 
Their interactions are non-trivial, as attested  
by Lemma \ref{Lemma-cancellation-2-Xi}.
To check \eqref{Check convergence-1}, 
we first note that
$$
A_{4,k}(t)=\sum_{j=1}^{3}A_{4,k,j}(t),\ \ \ A_{6,k}(t)=\sum_{i,j=1}^{3}\IP{A_{6,k,i}(t),A_{6,k,j}(t)}_{H^{\sigma}},$$
where 
\begin{equation*}
\left\{\begin{aligned}
&A_{4,k,1}(t)\triangleq 
\IP{(J^3_n-J^3_l)(a_k\A_k)^2 J_n u_n^{(R)}(t),\ v_{n,l}^{(R)}(t)}_{H^{\sigma}},\\
&A_{4,k,2}(t)\triangleq 
\IP{J^3_l (a_k\A_k)^2 (J_n-J_l) u_n^{(R)}(t),\ v_{n,l}^{(R)}(t)}_{H^{\sigma}},\\
&A_{4,k,3}(t)\triangleq 
\IP{J^3_l (a_k\A_k)^2 J_l v_{n,l}^{(R)}(t),\ v_{n,l}^{(R)}(t)}_{H^{\sigma}},\\
&A_{6,k,1}(t)\triangleq\, (J_n-J_l)(a_k\A_k)J_n u_n^{(R)}(t),\\
&A_{6,k,2}(t)\triangleq\, J_l (a_k\A_k) (J_n-J_l) u_n^{(R)}(t),\\
&A_{6,k,3}(t)\triangleq\, J_l a_k\A_k J_l v_{n,l}^{(R)}(t).
\end{aligned}\right.
\end{equation*}
Keep in mind that 
$\sigma\in\left(\frac{3}{2},{s-\max \{2\gamma_{0},1,2\theta{\bf 1}_{\{\eps>0\}}\}}\right)$ 
with $\gamma_0$ being given in  \eqref{gamma0}. 
Then we can infer from Hypothesis \ref{H-Q}, 
\eqref{OP boud Hs} and Lemma \ref{Lemma-Jn}  that
\begin{equation*}         
\sum_{k=1}^\infty A_{4,k,1}(t),\ 
\sum_{k=1}^{\infty}A_{4,k,2}(t)
\lesssim (l\land n)^{-(s-2\gamma_0-\sigma)} 
\big\|u_n^{(R)}(t)\big\|_{H^s}\big\|v_{n,l}^{(R)}(t)\big\|_{H^\sigma},              
\end{equation*}
\begin{equation*}
\sum_{k=1}^{\infty}\sum_{i,j\in \{1,2\}}
\IP{A_{6,k,i},A_{6,k,j}}_{H^{\sigma}} 
\lesssim  (l\land n)^{-2(s-\gamma_0-\sigma)} 
\big\|u_n^{(R)}(t)\big\|^2_{H^s},                   
\end{equation*}
\begin{equation*}
\sum_{k=1}^{\infty}\sum_{i\in \{1,2\}}
\IP{ A_{6,k,i},A_{6,k,3} }_{H^{\sigma}}
\lesssim
(l\land n)^{-(s-\gamma_0-\sigma)} 
\big\|u_n^{(R)}(t)\big\|_{H^s}\big\|v_{n,l}^{(R)}(t)\big\|_{H^\sigma}.      
\end{equation*}
Using Lemma \ref{Lemma-cancellation-2-Xi},
the remaining terms of $A_{4,k}$ and $A_{6,k}$ satisfies:
$$\sum_{k=1}^{\infty}\left\{A_{4,k,3}
+\IP{A_{6,k,3},A_{6,k,3}}_{H^{\sigma}}\right\}\lesssim \big\|v_{n,l}^{(R)}(t)\big\|_{H^\sigma}^2.$$
Consequently, we derive that for all $n,l\ge1$,
\begin{align*}
 \sum_{k=1}^\infty \left\{A_{4,k}(t)+A_{6,k}(t)\right\} 
\lesssim   \left(1+\big\|u_n^{(R)}(t)\big\|^2_{H^s}
+\big\|u_l^{(R)}(t)\big\|^2_{H^s}\right)
\left\{(l\land n)^{-2(s-2\gamma_0-\sigma)}
+\big\|v_{n,l}^{(R)}(t)\big\|_{H^\sigma}^2\right\}.
\end{align*}
Similarly, 
\begin{align*}
 \sum_{k=1}^\infty \left\{A_{5,k}(t)+A_{7,k}(t)\right\}  
\lesssim  \left(1+\big\|u_n^{(R)}(t)\big\|^2_{H^s}
+\big\|u_l^{(R)}(t)\big\|^2_{H^s}\right)
\left\{(l\land n)^{-2(s-2\gamma_0-\sigma)}
+\big\|v_{n,l}^{(R)}(t)\big\|_{H^\sigma}^2\right\}.
\end{align*}

For brevity, we omit the estimates for the other 
terms in \eqref{Check convergence-1} and 
also omit the proof for \eqref{Check convergence-2} 
since they can be handled similarly, 
and we refer to 
\cite{Tang-2023-JFA,
Tang-Wang-2022-arXiv,
Tang-Yang-2023-AIHP} for fuller exposition.

By  It\^{o}'s formula, the BDG inequality, 
\eqref{Check convergence-1}, and \eqref{Check convergence-2},  
we find a function 
$\varrho: \mathbb N\times\mathbb N \rightarrow[0,\infty)$ 
with $\displaystyle \lim_{n,l\to\infty} 
\varrho_{n,l}=0$ and constants $c_3>0$, $c_4=c_4(N)>0$ 
such that for all $n,l\ge 1$ and $t\in [0,N]$, 
\begin{align*} 
&\Ex\bigg[\sup_{t'\in [0,t\land \tau_N^{n,l}(R)]} \big\|v_{n,l}^{(R)}(t')\big\|_{H^\sigma}^2\bigg|\F_0\bigg]\\
\le \ & c_3 \Ex \bigg[\int_0^{t\land \tau_N^{n,l}(R)}\hh K(t',2N)\Big\{\varrho_{n,l}+\big\|v_{n,l}^{(R)}(t')\big\|_{H^\sigma}^2\Big\}\d t'\bigg|\F_0\bigg]\\
&+ c_3 \Ex\bigg[\bigg(\int_0^{t\land \tau_N^{n,l}(R)}\hh K(t',2N)\left\|v_{n,l}^{(R)}(t')\right\|_{H^\sigma}^2\Big\{\varrho_{n,l}
+\big\|v_{n,l}^{(R)}(t')\big\|_{H^\sigma}^2\Big\}\d t'\bigg)^{\frac 1 2}\bigg|\F_0\bigg]\\
\le \ &  \frac 1 2 \Ex\bigg[\sup_{t'\in [0,t\land \tau_N^{n,l}(R)]} \left\|v_{n,l}^{(R)}(t')\right\|_{H^\sigma}^2\bigg|\F_0\bigg]+ c_4 \varrho_{n,l} +c_4  \int_0^{t} \hh K(t',2N) \Ex\bigg[\sup_{r\in [0,t'\land\tau_N^{n,l}(R)]} \|v_{n,l}^{(R)}(r)\|_{H^\sigma}^2\bigg|\F_0\bigg]\d t'.
\end{align*}
By Gr\"onwall's inequality and noting $\varrho_{n,l}\to 0$ as $n,l\to\infty$, we arrive at \eqref{Cauchy in E-M}.

\ref{Convergence of u-n} We  fix a time $T>0$. 
On account of Lemma \ref{Lemma-un T estimates} and Chebyshev's inequality, 
we find that
for any $ N>T, n,l\ge 1$ and $\epsilon>0$,
\begin{align*}
&\p\bigg(\sup_{t\in[0,T]}\big\|v_{n,l}^{(R)}(t)\big\|_{H^\sigma}>\epsilon\bigg|\F_0\bigg)\\
\le \ &  \p(\tau_N^{n,l}(R)<T|\F_0) + \p\bigg(\sup_{t\in[0,\tau_N^{n,l}(R)]}\big\|v_{n,l}^{(R)}(t)\big\|_{H^\sigma}>\epsilon\bigg|\F_0\bigg)\\
\leq \ &\frac {2 Q(T,2R+\|u_0\|_{\Wlip})(1+\|u_0\|^2_{H^s})}{N^2}+\p\bigg(\sup_{t\in[0,\tau_N^{n,l}(R)]}\big\|v_{n,l}^{(R)}(t)\big\|_{H^\sigma}>\epsilon\bigg|\F_0\bigg).
\end{align*} 
According to  Lemma \ref{Lemma-un T estimates}, 
by first letting $n,l\to\infty$ and then $N\to\infty$, 
we obtain 
$$ \lim_{n,l\rightarrow\infty} \p\bigg(\sup_{t\in[0,T]}\big\|v_{n,l}^{(R)}(t)\big\|_{H^\sigma}>\epsilon\bigg|\F_0\bigg)=0,\ \ \epsilon,\,T>0.$$ 
By the Fatou lemma, this implies that for $\epsilon,\, T>0$,
\begin{align*} 
\limsup_{n,l\rightarrow\infty} \p\bigg(\sup_{t\in[0,T]}\big\|v_{n,l}^{(R)}(t)\big\|_{H^\sigma}>\epsilon\bigg)
\le  \Ex\bigg[\limsup_{n,l\rightarrow\infty} \p\bigg(\sup_{t\in[0,T]}\big\|v_{n,l}^{(R)}(t)\big\|_{H^\sigma}>\epsilon\bigg|\F_0\bigg)\bigg]=0.
\end{align*} 
Therefore, up to a subsequence, 
\eqref{Xn to X} holds for some progressively 
measurable process $u^{(R)}$. 
Together with 
Lemma \ref{Lemma-un T estimates}, this implies \eqref{X L2 bound}.
\end{proof}

\subsection{Step 3: concluding the proof for Theorem \ref{Thm:local-theory}}

Let $u^{(R)}(t)$ be given by 
\ref{Convergence of u-n} 
in Lemma \ref{Lemma-un convergence}.
For any $R\ge 1$, under the conditions 
in Theorem \ref{Thm:local-theory},  
we can use \ref{Convergence of u-n} in 
Lemma \ref{Lemma-un convergence} and 
Lemma  \ref{Lemma-un T estimates} to 
pass to the limit in \eqref{approximation cut eq} 
to obtain that for all $R\ge1$,
\begin{equation}\label{cut eq}
\left\{\begin{aligned}
\d u^{(R)}(t)&
+\Big[\eps \Lambda^{2\theta} u^{(R)}(t) 
+\lambda(t)\, u^{(R)}(t)
+G_{R}(u^{(R)}(t))\Big]\d t\\
&= \frac{1}{2}\sum _{k=1}^{\infty}
\left [ (a_k\A_{k})^{2} u^{(R)}(t) + (b_k\B_{k})^{2} u^{(R)}(t) \right ]
\d t+H_{R}\big(t,u^{(R)}(t)\big)\d \W(t),\\
u^{(R)}(0)&=u_0.
\end{aligned}\right.
\end{equation}
where
\begin{equation*}
G_{R}(u^{(R)}(t))\triangleq
\chi^2_R\big(\big\|u^{(R)}(t)-u_0\big\|_{\Wlip}\big)\Big[u^{(R)}(t)\,\pd_x u^{(R)}(t)+F(u^{(R)}(t))\Big]
\end{equation*} 
with $F(\cdot)$ being given in \eqref{F define}, 
and
\begin{equation*}
\left\{\begin{aligned}
&H_{R}\big(t,u^{(R)}(t)\big)e_{3k-2} 
\triangleq (a_k\A_k) u^{(R)}(t),\\
&H_{R}\big(t,u^{(R)}(t)\big)e_{3k-1} 
\triangleq (b_k\B_k) u^{(R)}(t),\\
&H_{R}\big(t,u^{(R)}(t)\big)e_{3k}
\triangleq \chi_R(\|u^{(R)}(t)-u_0\|_{\Wlip})h_k\big(t,u^{(R)}(t)\big).
\end{aligned}\right.
\end{equation*} 
Uniqueness of solution to \eqref{cut eq} follows 
from arguments leading to \eqref{Cauchy in E-M} 
and is omitted for brevity.

Now we are in a position to prove  Theorem \ref{Thm:local-theory}. 
Let 
\begin{equation*}
\tau^{(R)}\triangleq \inf\big\{t\ge 0: 
\big\|u^{(R)}(t)-u_0\big\|_{\Wlip}\ge R\big\}, 
\end{equation*}
and recall that 
$\sigma\in\left(\frac{3}{2},{s-\max \{2\gamma_{0},1,2\theta{\bf 1}_{\{\eps>0\}}\}}\right)$. 
By the continuity of $u^{(R)}(t)$ in 
$H^\sigma\hookrightarrow {\Wlip}$, 
we have $\p(\tau^{(R)}>0)=1$ for any $R>0$. 
Noting that $\chi_R(\|u^{(R)}(t)-u_0\|_{\Wlip})=1$ 
for $t\le \tau^{(R)}$, we see that \eqref{Cauchy problem-Ito} 
coincides with 
\eqref{cut eq} up to time $\tau^{(R)}$. 
This together with \eqref{X L2 bound} implies that 
$(u^{(R)},\tau^{(R)})$ is a local solution to 
\eqref{Cauchy problem-Ito}. 

Now we extend this local solution to a maximal solution as in Definition \ref{pathwise solution definition}.
By the uniqueness of solutions to \eqref{cut eq}, 
we see that $\tau^{(R)}$ is increasing in $R$, and 
$$u^{(R)}(t)= u^{(R+1)}(t),\ \ t\le \tau^{(R)}, \ \ R\ge1\ \ \p\text{-a.s.}$$
Let $\tau^* \triangleq\lim_{R\to\infty} \tau^{(R)}$, $\tau^{(0)}\triangleq0$ 
and we define
\begin{align*} 
u(t)\triangleq\sum_{R=1}^\infty \textbf{1}_{[\tau^{(R-1)}, \tau^{(R)})}(t) u^{(R)}(t),
\quad t\in [0,\tau^*).
\end{align*}
This allows us to conclude that $(u,\tau^*)$ is 
a local solution to \eqref{Cauchy problem-Ito}. 
Moreover, from the definitions of $\tau^*$ and $\tau^{(R)}$, 
\begin{equation*}
\limsup_{t\to\tau^*}\|u(t)\|_{\Wlip}
=\infty\ \text{on}\ \{\tau^*<\infty\} \quad \pas 
\end{equation*} 
Therefore $(u,\tau^*)$ is actually a maximal solution 
and \eqref{blow-up criterion-1} holds.

Finally, we demonstrate the continuity of the solution in time, 
i.e.,   $u \in C([0, \tau^*); H^s)$. 
As discussed in Remark \ref{Remark-singular}, 
we examine the mollified solution 
$J_n u$ ($J_n$ is defined in \eqref{Define Jn}).

Define the stopping time
\begin{equation}\label{TNN} 
\tau_N \triangleq N \land \inf \big\{ t \ge 0 : 
\| u(t) \|^2_{H^s} \ge N \big\}, \quad N \ge 1.
\end{equation}
By Lemma \ref{Lemma-un convergence}, 
it follows that $u \in C([0, \tau^*); H^\sigma)$. 
Given that $H^s \hookrightarrow H^\sigma$ is dense, 
$u$ is weakly continuous in $H^s$. 
Leveraging Hypotheses \ref{H-Q} and \ref{H-h} and applying Lemmas \ref{Lemma-cancellation-1}, 
\ref{Lemma-cancellation-2-Xi}, \ref{Lemma-Jn}, 
and \ref{uux+F Hs inner product Te}, we can identify 
a $\widetilde{K} \in \mathscr{K}$ (see \eqref{SCRK}) 
such that for all $N \ge 1$ and $t \ge 0$,
\begin{align*}
\sup_{n \ge 1,\, \| v \|_{H^s} \le N} 
\left| 2 \langle J_n [v \partial_x v + F(v) 
+ \lambda(t) v], J_n v \rangle_{H^s} 
+ \sum_{k=1}^{\infty} 
\| J_n h_k(t, v) \|_{H^s}^2 \right| \le \widetilde{K}(t, N),
\end{align*}
\begin{align*}
\sup_{n \ge 1,\, \| v \|_{H^s} \le N} 
\sum_{k=1}^{\infty} \Big(
\langle J_n a_k \mathcal{A}_k  v, J_n v \rangle_{H^s}^2 
+ \langle J_n b_k \mathcal{B}_k v, J_n v \rangle_{H^s}^2 
+ \langle J_n h_k(t, v), J_n v \rangle_{H^s}^2 
\Big) \le \widetilde{K}(t, N),
\end{align*}
and
\begin{align*}
\sup_{n \ge 1,\, \| v \|_{H^s} \le N} 
\bigg|\sum_{k=1}^{\infty}
\Big(\langle J_n (a_k \mathcal{A}_k)^2 v 
+ J_n (b_k \mathcal{B}_k)^2 v, J_n v \rangle_{H^s} 
+ \| J_n a_k \mathcal{A}_k v \|_{H^s}^2 
+ \| J_n b_k \mathcal{B}_k v \|_{H^s}^2 \Big) \bigg| 
\le \widetilde{K}(t, N).
\end{align*}
Applying It\^o's formula to $\| J_n u(t) \|_{H^s}^2$ and 
using the above estimates, for any $n \ge 1$, we find 
a constant $C_N > 0$ such that for some 
martingales $M^{(n)}$, the following estimates (understood in integral from) hold for $t \in [0, \tau_N]$ $\pas$,
\begin{equation}\label{Ito to Jn u-1}
 \mathrm{d} \langle M^{(n)} \rangle (t) 
\le C_N \, \mathrm{d} t,\quad
-C_N \, \mathrm{d} t \le \mathrm{d} \| J_n u(t) \|_{H^s}^2 
+ 2 \eps \| \Lambda^\theta J_n u(t) \|^2_{H^s} \, \mathrm{d} t
+ \mathrm{d} M^{(n)}(t) \le C_N \, \mathrm{d} t.
\end{equation}
Since $C_N$ in \eqref{Ito to Jn u-1} is $n$-independent, we can infer from the BDG inequality and
Fatou's lemma that
\begin{equation}\label{E X+GX F0}
\mathbb{E} \left[ \sup_{t \in [0, \tau_N]} \| u(t) \|^2_{H^s} 
+ 2 \eps\int_0^{\tau_N} \| \Lambda^\theta u(t) \|^2_{H^s} \, \mathrm{d} t 
\Bigg| \mathcal{F}_0 \right] < \infty, \quad N \ge 1.
\end{equation}
As a result, the stopping times
\begin{equation*}
\tau^\Gamma_N \triangleq N \land \inf \left\{ t \ge 0 : \Gamma(t) 
\triangleq \| u(t) \|^2_{H^s} 
+ 2\eps \int_0^{t} \| \Lambda^\theta u(t') \|^2_{H^s} \, \mathrm{d} t' \ge N 
\right\}, \quad N \ge 1,
\end{equation*}
satisfy $\tau^\Gamma_N \le \tau_N$ and
\begin{equation}\label{tau-Gamma-limit}
\mathbb{P} \left( \lim_{N \to \infty} \tau^\Gamma_N 
= \lim_{N \to \infty} \tau_N = \tau^* \right) = 1.
\end{equation}
Using \eqref{Ito to Jn u-1} and sending $n\to\infty$, we find that
\begin{equation*}
\mathbb{E} \left[ \Big| \Gamma(t \land \tau^\Gamma_N) 
- \Gamma(t' \land \tau^\Gamma_N) \Big|^4 \right] \le q(N) | t - t' |^2, 
\quad t, t' \ge 0, \  N \ge 1,
\end{equation*}
for some map $q : \mathbb{N} \to (0, \infty)$. 
By Kolmogorov's continuity criterion and \eqref{tau-Gamma-limit}, we see that  
$
\Gamma(\cdot) \in C([0, \tau^*))
$ $\pas$
Furthermore, \eqref{E X+GX F0} implies that 
the map $[0, \tau^*) \ni t \mapsto 
2 \int_0^t \| \Lambda^\theta u(t') \|^2_{H^s} \, \mathrm{d} t'$ 
is continuous $\pas$
Therefore, 
$\| u(\cdot) \|^2_{H^s} \in C([0, \tau^*))$ $\pas$ 
This and the weak continuity of $u$ in $H^s$ imply that $u \in C([0, \tau^*); H^s)$ $\pas$

\section{Diffusion-damping-noise interaction: Case 1}\label{Section:Diffusion-damping-noise-1}

The proof for \ref{T1:blow-up criterion} and \ref{T1:global}-\ref{T1:asymptotics} in Theorem \ref{Thm: global 1} will be presented in Sections \ref{Section:T1:blow-up-criterion}  and 
\ref{Section:T1:global-decay},  respectively. 

\subsection{Refined blow-up criterion}
\label{Section:T1:blow-up-criterion}

\begin{proof}[Proof for \ref{T1:blow-up criterion}  
in Theorem \ref{Thm: global 1}]

The proof is divided into two steps: 
\begin{itemize}
	\item refine \eqref{blow-up criterion-1} by replacing $\limsup_{t\rightarrow \tau^*}\|u(t)\|_{\Wlip}
	=\infty$ with $\limsup_{t\rightarrow \tau^*}\int_0^t\|u(t')\|_{\Wlip}\d t'
	=\infty$ (see \eqref{blow-up criterion-3} below);
	
	\item demonstrate that  $\limsup_{t\rightarrow \tau^*}\int_0^t\|u(t')\|_{L^\infty}\d t'
	<\infty$
\end{itemize}

\textbf{Step 1.}  We claim that the following 
blow-up criterion holds:
\begin{equation}\label{blow-up criterion-3}
\begin{aligned}
\mathbf{1}_{\left\{
\limsup_{t \rightarrow \tau^*} \| u(t) \|_{H^s} = \infty \right\}}
= \mathbf{1}_{\left\{ \limsup_{t \rightarrow \tau^*} \int_{0}^{t}
\| u(t') \|_{\Wlip} \, \mathrm{d} t' = \infty \right\}} \
\text{on} \ \{ \tau^* < \infty \} \ \mathrm{a.s.}
\end{aligned}
\end{equation}
Recall that 
$s > 3/2 + \max \{ 2 \gamma_0, 1, 2 \theta \mathbf{1}_{\{ \eps > 0 \}} \}$.  Let $\eta \in \left[ 1, s - 
\max \{ 2 \gamma_0, 1, 2 \theta \mathbf{1}_{\{ \eps > 0 \}} \} \right)$.
Using It\^o's formula,
we have:
\begin{align}
\mathrm{d} \| u(t) \|^2_{H^\eta} 
= \ &   \left\{-2\eps \norm{\Lambda^{\theta}u(t)}^2_{H^\eta}-2\lambda(t)\norm{u(t)}^2_{H^\eta}
-2\bIP{u(t) \partial_x u(t) 
+ F(u(t)), u(t)}_{H^\eta}\right\} \mathrm{d} t
+ \sum_{i=1}^{5} I_{i,\eta}(t) \, \mathrm{d} t \notag \\
& + 2 \sum_{k=1}^{\infty} E_{1,\eta,k}(t) \, \mathrm{d} \overline{W}_k(t)
+ 2 \sum_{k=1}^{\infty} E_{2,\eta,k}(t) \, \mathrm{d} \widehat{W}_k(t)
+ 2 \sum_{k=1}^{\infty} E_{3,\eta,k}(t) \, \mathrm{d} \widetilde{W}_k(t), 
\label{Ito-u-Hsigma}
\end{align} 
where
\begin{align}\label{I+E eta}
\begin{cases} 
I_{1,\eta}(t)\triangleq\sum_{k=1}^{\infty}
\IP{(a_k\A_k)^2 u(t), u(t)}_{H^{\eta}},
\ \
&I_{2,\eta}(t)\triangleq  \sum_{k=1}^{\infty}
\norm{a_k\A_ku(t)}^2_{H^{\eta}},\\[10pt]
I_{3,\eta}(t)\triangleq\sum_{k=1}^{\infty}
\IP{(b_k\B_k)^2 u(t), u(t)}_{H^{\eta}},
\ \
&I_{4,\eta}(t)\triangleq  \sum_{k=1}^{\infty}
\norm{b_k\B_ku(t)}^2_{H^{\eta}},\\[10pt]
I_{5,\eta}(t)\triangleq\sum_{k=1}^{\infty}
\norm{h_k(t,u(t))}^2_{H^{\eta}},  \ \ 
&E_{1,\eta,k}(t)\triangleq\bIP{a_k\A_k  u(t), u(t)}_{H^{\eta}},\\[10pt]
E_{2,\eta,k}(t)\triangleq\bIP{b_k\B_k  u(t), u(t)}_{H^{\eta}},\ \ 
&E_{3,\eta,k}(t)\triangleq\bIP{h_k(t,u(t)), u(t)}_{H^{\eta}}.
\end{cases} 
\end{align}
Using Hypotheses \ref{H-Q} and \ref{H-h} with 
$K_1 \equiv c_0>0$ and applying Lemma \ref{Lemma-cancellation-2-Xi}, 
we derive the following estimates:
\begin{equation}\label{eq:estm_A1A2HqHh}
\left\{ \begin{aligned}
& \left[ \left| I_{1,\eta} 
+ I_{2,\eta} \right| 
+ \left| I_{3,\eta} 
+ I_{4,\eta} \right| \right](t) \lesssim \| u(t) \|^2_{H^\eta},  \\
&\left| I_{5,\eta} \right|(t) \leq c_0\left(1 + \| u(t) \|^2_{H^\eta} \right), \\ 
&\sum_{k=1}^{\infty} \left[ E_{1,\eta,k}^2 
+ E_{2,\eta,k}^2 \right](t) \lesssim \| u(t) \|_{H^\eta}^4, \\
& \sum_{k=1}^{\infty} E_{3,\eta,k}^2(t) 
\leq c_0\left(1 + \| u(t) \|^2_{H^\eta} \right) \| u(t) \|^2_{H^\eta}.
\end{aligned} \right.
\end{equation}

Let $\sigma \in \left( \frac{3}{2}, 
s - \max \{ 2 \gamma_0, 1, 2 \theta \mathbf{1}_{\{ \eps > 0 \}} \} \right)$. 
Using Lemmas \ref{F Lemma} and \ref{Taylor-commutator}, 
we find that there exists a constant $C_0 > 0$ such that
\begin{align*}
\left| 2 \langle u(t) \partial_x u(t) 
+ F(u(t)), u(t) \rangle_{H^\sigma} \right|
\leq C_0 \| u(t) \|_{\Wlip} \| u(t) \|^2_{H^\sigma}.
\end{align*}

Next, for $\sigma\in\left(\frac{3}{2},{s-\max \{2\gamma_{0},1,2\theta{\bf 1}_{\{\eps>0\}}\}}\right)$,  applying the It\^o formula to $\log({\rm e} + \| u \|^2_{H^\sigma})$, 
neglecting the damping and diffusion terms (which have the appropriate signs), alongside the negative It\^o correction term,  and
substituting the estimates in \eqref{eq:estm_A1A2HqHh} 
with $\eta = \sigma$, we obtain that for some constant $C_1 > 0$,
\begin{align}
 \mathrm{d} \log({\rm e} + \| u(t) \|^2_{H^\sigma}) 
\leq \ &   C_1 \| u(t) \|_{\Wlip} \frac{\| u(t) \|^2_{H^\sigma}}{ {\rm e} 
+ \| u(t) \|^2_{H^\sigma} } \, \mathrm{d} t 
+ C_1 \frac{ 1 + \| u(t) \|^2_{H^\sigma} }{ {\rm e} 
+ \| u(t) \|^2_{H^\sigma} } \, \mathrm{d} t \notag \\
& + 2 \frac{ \sum_{k=1}^{\infty} 
\left( E_{1,\sigma,k}(t) \, \mathrm{d} \overline{W}_k(t) 
+ E_{2,\sigma,k}(t) \, \mathrm{d} \widehat{W}_k(t) 
+ E_{3,\sigma,k}(t) \, \mathrm{d} \widetilde{W}_k(t) \right) }
{ {\rm e} + \| u(t) \|^2_{H^\sigma} }.\label{Ito-log(u)}
\end{align}
Define
\begin{equation*}
\widetilde{\tau} \triangleq \lim_{N \to \infty} \widetilde{\tau}_N, 
\quad 
\widetilde{\tau}_N \triangleq N \land 
\inf \left\{ t \ge 0 : \int_0^t \| u(t') \|_{\Wlip} \, \mathrm{d} t' 
\ge N \right\}, 
\quad N \ge 1,
\end{equation*}
and
\begin{equation*}
\widehat{\tau} \triangleq \lim_{M \to \infty} \widehat{\tau}_M, 
\quad 
\widehat{\tau}_M \triangleq M \land 
\inf \left\{ t \ge 0 : \| u(t) \|_{H^\sigma} \ge M \right\}, 
\quad M \ge 1.
\end{equation*}
Let $T>0$.
Using the BDG inequality, 
\eqref{eq:estm_A1A2HqHh}, 
we can find 
constants $C_2,C_3>0$ such that
\begin{align*}
&\Ex\left[\sup_{t\in[0,T\land\widetilde{\tau}_N\land\hh\tau_M]}
\log\left({\rm e}+\norm{u(t)}^2_{H^{\sigma}}\right)  
\bigg|\F_0\right]\\
\leq  \ &\log({\rm e}+\norm{u_0}^2_{H^{\sigma}})
+C_2N+C_2T +C_2 \Ex\Bigg[
\Bigg(\int_{0}^{T\land \widetilde{\tau}_N\land\hh\tau_M}
\frac{\big(1+\norm{u(t)}^2_{H^\sigma}\big)^2}{\big({\rm e}+\norm{u(t)}^2_{H^{\sigma}}\big)^2}
\d  t\Bigg)^{\frac12}
\Bigg|\F_0\Bigg]\\
\leq  \ &\log({\rm e}+\norm{u_0}^2_{H^{\sigma}}) 
+C_3N +C_3(\sqrt{T}+T),\quad N,M\ge1.
\end{align*} 
Therefore, for any \( M, N \ge 1 \) and \( T > 0 \),
\begin{align*}
\mathbb{P}\Big(\widehat{\tau} 
\le T \land \widetilde{\tau}_N \big| \mathcal{F}_0\Big)
\leq \mathbb{P}\Big(\widehat{\tau}_M 
\le T   \land \widetilde{\tau}_N \big| \mathcal{F}_0\Big) 
\leq \frac{ \log({\rm e} + \| u_0 \|^2_{H^\sigma}) 
+ C_3 N + C_3 (\sqrt{T} + T) }{\log({\rm e} + M^2)}.
\end{align*}
Sending \( M \to \infty \), we see that 
\[ \mathbb{P}\Big(\widehat{\tau} > 
T \land \widetilde{\tau}_N \big| \mathcal{F}_0\Big) = 1 \]
and therefore 
\[
\mathbb{P}\Big(\widehat{\tau} > T \land \widetilde{\tau}_N\Big) 
= \mathbb{E}\Big[\mathbb{P}
\Big(\widehat{\tau} > T \land \widetilde{\tau}_N \big| 
\mathcal{F}_0\Big)\Big] 
= 1, \quad N \ge 1, \, T > 0.
\]
By the arbitrariness of $T > 0$, we then derive:
\[
\mathbb{P}\Big(\widehat{\tau} \ge \widetilde{\tau}\Big) 
= \mathbb{P}\bigg( \bigcap_{k \ge 1} \bigcap_{N \ge 1} 
\{\widehat{\tau} > k \land \widetilde{\tau}_N\} \bigg) = 1.
\]
Conversely, 
since \( H^\sigma \hookrightarrow \Wlip \), 
it is easy to see that for some \( C >0 \),
\[
\int_0^{\widehat{\tau}_M} \| u(t') \|_{\Wlip} \, \mathrm{d} t' 
\leq C \int_0^M M \, \mathrm{d} t' \leq \mathrm{ceil}(CM^2) \triangleq M_0,
\]
where \(\mathrm{ceil}(x)\) denotes the ceiling function of \( x \in \mathbb{R} \).
As a result, we obtain 
$\mathbb{P}(\widehat{\tau}_M 
\le \widetilde{\tau}_{M_0} \le \widetilde{\tau}) = 1$ 
for all  $M \ge 1$ 
and hence
\[
\mathbb{P}(\widehat{\tau} \le \widetilde{\tau}) = 1.
\]
In conclusion, we have 
\( \mathbb{P}(\widehat{\tau} = \widetilde{\tau}) = 1 \), 
which, together with \eqref{blow-up criterion-1}, 
implies  \eqref{blow-up criterion-3}.

\textbf{Step 2.} 
With \eqref{blow-up criterion-3} in hand, we proceed to prove \eqref{blow-up criterion-2} 
by demonstrating that 
\begin{equation}\label{to show bc2}
\limsup_{t \to \tau^*} \int_0^t \| u(t') \|_{L^\infty} \, \mathrm{d} t' 
< \infty \quad \text{on} \ \{ \tau^* < \infty \}. 
\end{equation}

Recall \eqref{TNN}, i.e.,
\begin{equation*}
\tau_N \triangleq N \land \inf \big\{ t \ge 0 : \| u(t) \|^2_{H^s} \ge N \big\}, 
\quad N \ge 1.
\end{equation*}
Let \( T > 0 \). We can infer from \eqref{Ito-u-Hsigma} 
with \( \eta = 1 \) and Lemma \ref{H1-conserve} that there 
is a constant \( C_4 > 0 \) such that
\begin{align}
&{\rm d}\norm{u(t)}_{H^1}^2 
+2\eps \norm{\Lambda^\theta u(t)}_{H^1}^2\d t
+2 \lambda(t)\norm{u(t)}_{H^1}^2\d t \notag\\
\leq \ &
C_4\left(1+\norm{u(t)}^2_{H^1}\right)\d t 
+2\sum_{k=1}^{\infty}E_{1,1,k}(t)\d  \ol W_k(t)  
 +2\sum_{k=1}^{\infty}E_{2,1,k}(t)\d  \hh W_k(t)
+2\sum_{k=1}^{\infty}E_{3,1,k}(t)\d  \widetilde{W}_k(t).\label{Ito-H1}
\end{align}
Applying the BDG inequality on 
\( t \in [0, T \land \tau_N] \) 
for \( N \ge 1 \), and using \eqref{eq:estm_A1A2HqHh}, 
we find a constant \( C_5 > 0 \) such that
\begin{align*}
& \mathbb{E} \left[ \sup_{t \in [0, T \land \tau_N]} 
\| u(t) \|_{H^1}^2 \Big| \mathcal{F}_0 \right] 
+ 2 \eps\mathbb{E} 
\left[ \int_0^{T \land \tau_N} 
\| \Lambda^\theta u(t) \|_{H^1}^2 \, \mathrm{d} t \bigg| \mathcal{F}_0 \right] \notag\\
\leq  \ &\| u_0 \|^2_{H^1} + \frac{1}{2} \mathbb{E} 
\left[ \sup_{t \in [0, T \land \tau_N]} \| u(t) \|_{H^1}^2 \Big| \mathcal{F}_0 \right] + C_5 \int_0^T \mathbb{E} 
\left[ 
\left(1 + \sup_{t' \in [0, t \land \tau_N]} \| u(t') \|_{H^1}^2 \right) 
\bigg| \mathcal{F}_0 \right] \, \mathrm{d} t, \quad N \ge 1.
\end{align*}
On account of  Gr\"onwall's inequality, there exists a function 
\( \bar{K} : [0, \infty) \times [0, \infty) \to (0, \infty) \) 
that is increasing in both variables such that 
for every \( N \ge 1 \),
\begin{align*}
\mathbb{E} \left[ \sup_{t \in [0, T \land \tau_N]} \| u(t) \|_{H^1}^2 
+ 4 \eps \int_0^{T \land \tau_N} \| \Lambda^\theta u(t) \|_{H^1}^2 
\, \mathrm{d} t \bigg|\F_0 \right]
\leq \bar{K}(T, \| u_0 \|_{H^1}).
\end{align*}
By letting \( N \to \infty \), applying  Beppo Levi's lemma, 
and utilizing the time continuity of \( u \in H^1 \), we have
\begin{align*}
\mathbb{E} \left[ \sup_{t \in [0, T \land \tau^*]} \| u(t) \|_{H^1}^2 
+ 4 \eps \int_0^{T \land \tau^*} \| \Lambda^\theta u(t') \|_{H^1}^2 
\, \mathrm{d} t' \bigg|\F_0 \right]
\leq \bar{K}(T, \| u_0 \|_{H^1}).
\end{align*}
Then \eqref{to show bc2} follows
from the above estimate and the embedding \( H^1 \hookrightarrow L^\infty \). 
Combining the two steps, we  
obtain \eqref{blow-up criterion-2}.
\end{proof}

\begin{remark}
Two remarks concerning the applicability of Gr\"onwall's inequality in the above proofs are as follows:

\begin{enumerate}
\item From the context, \eqref{Ito to Jn u-1}, \eqref{Ito-u-Hsigma}, \eqref{Ito-log(u)} and \eqref{Ito-H1} should be understood in integral form, since the damping and diffusion terms are neglected and the integral form of Gr\"onwall's inequality applies.
\item In the following, when the damping effect requires a more precise analysis, the positivity condition in the integral form of Gr\"onwall's inequality may not be satisfied. In that case, we use the differential form of Gr\"onwall's inequality; see \eqref{eq:gronwall_differntial} below.
\end{enumerate}
\end{remark}

\subsection{Global existence \&  asymptotics}
\label{Section:T1:global-decay}

We are now in a position to prove global existence and asymptotic behaviour for solutions under the linear growth noise condition.

\begin{proof}[Proof for  \ref{T1:global} and \ref{T1:asymptotics}
in Theorem \ref{Thm: global 1}] 
Recall the stopping time $\tau_N$ defined in \eqref{TNN}:
\begin{equation*}
\tau_N \triangleq N \land \inf \big\{ t \ge 0 : \| u(t) \|^2_{H^s} \ge N \big\}, \quad N \ge 1,
\end{equation*}
and consider the function
\begin{equation*} 
G(\bullet\land \tau_N) \triangleq \mathbb{E}\left[\|u(\bullet\land \tau_N)\|_{H^1}^2 + 2\eps \int_0^{\bullet\land \tau_N} \|\Lambda^\theta u(t')\|_{H^1}^2 \mathrm{d} t' 
+ 2\int_0^{\bullet\land \tau_N} \lambda(t') \|u(t')\|_{H^1}^2 \mathrm{d} t' \bigg|\mathcal{F}_0\right],\quad N\ge1.
\end{equation*}

Based on \eqref{Ito-u-Hsigma} (understood in integral form), \eqref{I+E eta} with $\eta=1$, Hypothesis \ref{H-h} with $K_1 \equiv c_0 > 0$, and Lemmas \ref{H1-conserve} and \ref{Lemma-cancellation-2-Xi}, we obtain for all $t_1 > t_2 \ge 0$:
\begin{align*}
&|G(t_1 \land \tau_N) - G(t_2 \land \tau_N)| \\
\leq\ & \mathbb{E}\left[\int_{t_2}^{t_1} \mathbf{1}_{\{t' \le \tau_N\}} \big[ |I_{1,1} + I_{2,1}| + |I_{3,1} + I_{4,1}| + |I_{5,1}| \big](t') \mathrm{d} t' \bigg|\mathcal{F}_0\right] \\
\leq\ &  (\Xi + c_0) \mathbb{E}\left[\int_{t_2}^{t_1} \mathbf{1}_{\{t' \le \tau_N\}} \left(1 + \|u(t')\|_{H^1}^2\right) \mathrm{d} t' \bigg|\mathcal{F}_0\right] \\
\leq\ & (\Xi + c_0)(1 + N^2)(t_1 - t_2).
\end{align*}

This implies that the mapping $t \mapsto G(t \land \tau_N)$ is locally Lipschitz, and hence absolutely continuous on $[0,\infty)$. 
Since 
\begin{equation*}
\mathbb{E}\left[2\eps \int_0^{\bullet\land \tau_N} \|\Lambda^\theta u(t')\|_{H^1}^2 \mathrm{d} t' 
+ 2\int_0^{\bullet\land \tau_N} \lambda(t') \|u(t')\|_{H^1}^2 \mathrm{d} t' \bigg|\mathcal{F}_0\right]\in AC_{\rm loc}([0,\infty)) \quad \pas,
\end{equation*}
we conclude that
\begin{equation}\label{Ex u differential}
\mathbb{E}\left[\|u(\bullet\land \tau_N)\|_{H^1}^2 \Big|\mathcal{F}_0\right] \quad \text{is absolutely continuous} \quad \pas
\end{equation}

Consequently, from \eqref{Ito-u-Hsigma} (understood in integral form),  \eqref{I+E eta} with $\eta=1$, \eqref{Ex u differential}, Hypothesis \ref{H-h} with $K_1 \equiv c_0 > 0$, and Lemmas \ref{H1-conserve} and \ref{Lemma-cancellation-2-Xi}, we derive the differential inequality:
\begin{align}
&\frac{\mathrm{d}}{\mathrm{d}t} \mathbb{E}\left[\|u(t\land \tau_N)\|_{H^1}^2 + 2\eps \int_0^{t\land \tau_N} \|\Lambda^\theta u(t')\|_{H^1}^2 \mathrm{d} t'
\Big|\mathcal{F}_0\right]\notag\\
\leq\ & c_0 +  (\Xi + c_0) \mathbb{E}\left[\|u(t\land \tau_N)\|_{H^1}^2 \Big|\mathcal{F}_0\right]-
2\mathbb{E}\left[\mathbf{1}_{\{t\le \tau_N\}} \lambda(t)\|u(t)\|_{H^1}^2 \Big|\mathcal{F}_0\right],\quad t\ge0, \quad N\ge1,\label{eq:gronwall_differntial}
\end{align}
where $\Xi$ is given in Lemma \ref{Lemma-cancellation-2-Xi}.

We now prove $\mathbb{P}(\tau^* = \infty) = 1$. First, neglecting the damping term in \eqref{eq:gronwall_differntial} yields
\begin{align*}
 \frac{\mathrm{d}}{\mathrm{d}t}\mathbb{E}\left[\|u(t\land \tau_N)\|_{H^1}^2 + 2\eps \int_0^{t\land \tau_N} \|\Lambda^\theta u(t')\|_{H^1}^2 \mathrm{d} t' \bigg|\mathcal{F}_0\right] 
 \leq  c_0 + (\Xi + c_0)\mathbb{E}\left[\|u(t\land \tau_N)\|_{H^1}^2 + 2\epsilon \int_0^{t\land \tau_N} \|\Lambda^\theta u(t')\|_{H^1}^2 \mathrm{d} t' \bigg|\mathcal{F}_0\right],
\end{align*}
which implies the bound
\begin{align*}
\mathbb{E}\left[\|u(t\land \tau_N)\|_{H^1}^2 + 2\eps \int_0^{t\land \tau_N}\|\Lambda^\theta u(t')\|_{H^1}^2 \mathrm{d} t' \bigg|\mathcal{F}_0\right]  \leq \|u_0\|_{H^1}^2 e^{(\Xi + c_0)t} + \frac{c_0}{\Xi + c_0}\left(e^{(\Xi + c_0)t} - 1\right).
\end{align*}
Since $\lim_{N\to\infty}\tau_N=\tau^*$ $\pas$, 
from the above estimate and Lemma \ref{Fractional Laplacian}, and applying the Cauchy-Schwarz inequality in time, we obtain
\[
\limsup_{t \rightarrow \tau^*}\int_0^{t} \|\partial_x u(t')\|_{L^\infty} \mathrm{d} t' < \infty  \ {\rm on}\ \{\tau^*<\infty\},
\]
which, together with \eqref{blow-up criterion-2}, proves $\mathbb{P}(\tau^* < \infty) = 0$.

 When including the damping effect, we first obtain an exact integral identity for $\mathbb{E}[\|u(t\land \tau_N)\|_{H^1}^2 |\mathcal{F}_0]$ from \eqref{Ito-u-Hsigma} with $\eta=1$. Next, we send $N\to\infty$ (hence $\tau_N\to\tau^*=\infty$ $\pas$). We then differentiate the resulting absolutely continuous identity and apply the bounds in \eqref{eq:estm_A1A2HqHh} to find
 \begin{align*}
 	\frac{\mathrm{d}}{\mathrm{d}t} \mathbb{E}\left[\|u(t)\|_{H^1}^2 + 2\eps \int_0^{t} \|\Lambda^\theta u(t')\|_{H^1}^2 \mathrm{d} t'
 	\Big|\mathcal{F}_0\right] 
 	\leq   c_0 + \big((\Xi + c_0)-2\lambda(t) \big) \mathbb{E}\left[\|u(t)\|_{H^1}^2 \Big|\mathcal{F}_0\right],\quad t\ge0.
 \end{align*}
 Let $B(t)\triangleq(\Xi+c_0) t -2 \int_0^t\lambda(r) \d r$. Applying Gr\"onwall's inequality in differential form to the above inequality leads to
 \begin{align*}
 	\Ex\left[\norm{u(t)}_{H^1}^2\Big|\F_0\right] 
 	\leq   \norm{u_0}_{H^1}^2
 	\exp \big(B(t)\big) +c_0\int_0^t\exp \big(B(t)-B(t')\big)\d t'.
 \end{align*} 
 The two cases \eqref{limsup Ex u} and \eqref{lim Ex u=0} follow directly from this inequality, completing the proof.
\end{proof}

\section{Diffusion-damping-noise interaction: Case 2}\label{Section:Diffusion-damping-noise-2}
In this section we prove Theorem \ref{Thm: global 2}. 
The proofs for \ref{T2:global}, \ref{T2:stability}, 
and \ref{T2:invariant measure} in 
Theorem \ref{Thm: global 2} are given in 
Sections \ref{Section:T2:global}, \ref{Section:T2:stability}, 
and  \ref{Section:T2:invariant measures}, respectively.

\subsection{Global solution}\label{Section:T2:global}

We will now prove the second global existence result using a Lyapunov estimate on $u(t)$.

\begin{proof}[Proof for \ref{T2:global} in Theorem \ref{Thm: global 2}]

Recall $J_n$ in \eqref{Define Jn}.  
As mentioned in Remark \ref{Remark-singular} 
and similar to \eqref{Ito-u-Hsigma},   
we apply It\^o's formula to $\norm{J_n u(t)}^2_{H^{s}}$ 
to find 
\begin{align}
&\norm{J_n u(t)}^2_{H^{s}}-\norm{J_n u_0}^2_{H^{s}}
+\int_0^t \left\{2\eps\norm{J_n \Lambda^\theta  u(t)}^2_{H^s}
+2\lambda(t)\norm{J_n u(t)}^2_{H^s}\right\}\d t\notag\\
=\ & \int_0^t \left\{- 2\bIP{J_n [u(t)\,\pd_xu(t)]+ J_n F(u(t)),
J_n u(t)}_{H^{s}}+\sum_{i=1}^{5}I_{i,s}^{(n)}(t)\right\} \d t 
\notag\\
&+2\int_0^t\sum_{k=1}^{\infty}E_{1,s,k}^{(n)}(t)\d  \ol W_k(t) 
+2\int_0^t\sum_{k=1}^{\infty}E_{2,s,k}^{(n)}(t)\d  \hh W_k(t)
+2\int_0^t\sum_{k=1}^{\infty}E_{3,s,k}^{(n)}(t)\d  \widetilde{W}_k(t),
\label{Ito to Jn u-2}
\end{align} 
where
\begin{align*}
\begin{cases} 
I_{1,s}^{(n)}(t)\triangleq\sum_{k=1}^{\infty}
\IP{J_n(a_k\A_k)^2 u(t), J_nu(t)}_{H^{s}},
\ \
&I_{2,s}^{(n)}(t)\triangleq\sum_{k=1}^{\infty}
\norm{J_na_k\A_ku(t)}^2_{H^{s}},\\[10pt]
I_{3,s}^{(n)}(t)\triangleq\sum_{k=1}^{\infty}
\IP{J_n(b_k\B_k)^2 u(t), J_nu(t)}_{H^{s}},
\ \
&I_{4,s}^{(n)}(t)\triangleq\sum_{k=1}^{\infty}
\norm{J_n b_k\B_ku(t)}^2_{H^{s}},\\[10pt]
I_{5,s}^{(n)}(t)\triangleq\sum_{k=1}^{\infty}
\norm{J_n h_k(t,u(t))}^2_{H^{s}},  \ \
&E_{1,s,k}^{(n)}(t)
\triangleq\bIP{J_n a_k\A_k  u(t), J_n u(t)}_{H^{s}},\\[10pt]
E_{2,s,k}^{(n)}(t)
\triangleq\bIP{J_n b_k\B_k  u(t), J_n u(t)}_{H^{s}},\ \ 
&E_{3,s,k}^{(n)}(t)
\triangleq\bIP{J_n h_k(t,u(t)), J_n u(t)}_{H^{s}}.
\end{cases} 
\end{align*}
Similar to \eqref{eq:estm_A1A2HqHh}, 
we infer from Hypothesis \ref{H-Q}, 
Lemmas \ref{Lemma-Jn} 
and \ref{Lemma-cancellation-2-Xi} 
that  
\begin{align*}
\Big[ \abs{I_{1,s}^{(n)}+I_{2,s}^{(n)}}+\abs{I_{3,s}^{(n)}
+I_{4,s}^{(n)}}\Big](t)\leq \Xi  \norm{u(t)}^2_{H^s},\quad  \sum_{k=1}^{\infty}\left[\left(E_{1,s,k}^{(n)}\right)^2
+\left(E_{2,s,k}^{(n)}\right)^2\right](t)
\lesssim \norm{u(t)}_{H^s}^4.
\end{align*}
Recall the stopping time $\tau_N$ in \eqref{TNN}, i.e.,
\begin{equation*} 
\tau_N\triangleq
N\land \inf\big\{t\ge 0: \|u(t)\|^2_{H^{s}}\ge N\big\},
\quad N\ge 1.
\end{equation*}

Applying It\^o's formula to $V(\|J_n u(t)\|^2_{H^s})$ and discarding the non-positive It\^o correction terms involving  $V'' \le 0$, we  then use 
Hypotheses \ref{H-h} and \ref{H-h-large},  Lemmas 
\ref{uux+F Hs inner product Te} 
and \ref{Lemma-cancellation-2-Xi} to derive 
that for any $t>0$,
\begin{align*}
&\Ex \left[V(\|J_n u(t\wedge \tau_N)\|^2_{H^s})\big|\F_0\right]
-V(\|u_0\|^2_{H^s})  \\
\leq \ &\Ex\bigg[
\int_0^{t\wedge \tau_N}V'(\|J_n u(t')\|^2_{H^s})
\Big\{
\left(\Xi+2\Theta \|u(t')\|_{W^{1,\infty}}\right)
\|u(t')\|^2_{H^s} \Big\}\, {\rm d}t'\Big|\F_0\bigg]\\
&- \Ex\bigg[
\int_0^{t\wedge \tau_N} 2\lambda(t')V'(\|J_n u(t')\|^2_{H^s})
\|J_n u(t')\|^2_{H^s}\, {\rm d}t'\Big|\F_0\bigg]\\
&+ \Ex\bigg[
\int_0^{t\wedge \tau_N}V'(\|J_n u(t')\|^2_{H^s})
I_{5,s}^{(n)}(t')\, {\rm d}t'\Big|\F_0\bigg]\\
& +\Ex\left[\int_0^{t\wedge \tau_N}
2V''(\|J_n u(t')\|^2_{H^s})\left(\sum_{k=1}^\infty 
\left(E_{3,s,k}^{(n)}(t')\right)^2 \right)\, {\rm d}t'\Big|\F_0\right].
\end{align*}
Letting $n \to \infty$, by Fatou's lemma,  the dominated convergence theorem 
and Hypothesis \ref{H-h-large}, one has
\begin{align*}
  \Ex \left[V(\|u(t\wedge \tau_N)\|^2_{H^s})\big|\F_0\right]
-V(\| u_0\|^2_{H^s})  
\leq   
\int_0^{t}g(t')\Ex\left[V(\|u(t'\land\tau_N)\|^2_{H^s})
\Big|\F_0\right]\, {\rm d}t',
\end{align*}
which implies
\begin{equation}\label{Hs bound V tau-N}
\sup_{N\ge 1}\Ex \left[V(\|u(t\wedge \tau_N)\|^2_{H^s})\Big|\F_0\right]
\leq V(\| u_0\|^2_{H^s}) {\rm e}^{\int_0^t g(t')\d t'}, 
\quad t>0.
\end{equation}
Accordingly, leveraging the time continuity of $u(t)$ in $H^s$, 
we deduce that for $N\ge t>0$,
\begin{align*} 
\p\Big(\tau^*<t\big|\F_0\Big)  \le\,  \p\Big(\tau_{N}<t\big|\F_0\Big) 
\le \frac{ V(\|u_0\|_{H^s}^{2}){\rm e}^{\int_0^tg(t')\d t'}}{V(N)}.
\end{align*} 
Letting $N\uparrow\infty$ and then $t\uparrow\infty$ 
yields
$\p(\tau^*<\infty|\F_0)=0$. 
Hence we obtain $\p(\tau^*<\infty)=0.$ Then \eqref{Hs bound V} is a consequence of \eqref{Hs bound V tau-N}.
\end{proof}

\subsection{Continuity of the solution map}
\label{Section:T2:stability}

Recall the definition \eqref{X-s1s2} of the spaces $X_{s_1,s_2}$.
With \eqref{Hs bound V}, we can prove that the solution map is continuous from $X_{s,\sigma}$ to $X_{s,\sigma}$.

\begin{proof}[Proof for \ref{T2:stability} in Theorem \ref{Thm: global 2}]

Define the stopping times:
\begin{equation*}
\tau_N^{n}\triangleq\inf\{t\ge 0: \|u_{n}(t)\|_{H^s}
\lor\|u(t)\|_{H^s}\ge N\},\ \ n, N\ge1.
\end{equation*}

We note that $u_n$ and $u$ are continuous in $t$. 
Therefore, by Chebyshev's inequality, $\tau_N^n$ 
satisfies the following for every $n,N\ge 1$:
\begin{align*}
\p\left(\tau_N^n<t\big|\F_0\right) \leq 
\frac{\Big(V(\|u_{0,n}\|^2_{H^s})+V(\|u_0\|^2_{H^s})\Big){\rm e}^{\int_0^tg(t')\d t'}}{V(N^2)}.
\end{align*}
Using Lemma \ref{Taylor-commutator} leads 
to the following
inequality:
\begin{align*}
\left|\bIP{u\,\pd_xu-v\,\pd_x v,u-v}_{H^{\sigma}}\right|
\lesssim \left(\|u\|_{H^{\sigma+1}}+\|v\|_{H^{\sigma+1}}\right)
\|u-v\|^2_{H^\sigma},\ \ u,v\in H^{s}.
\end{align*}
Combining this,  Lemmas \ref{F Lemma} 
and \ref{Lemma-cancellation-2-Xi} 
and Hypothesis  \ref{H-h}, we find a constant $C>0$ 
such that for all $t>0$ and $n,N\ge1$,
\begin{align*} 
&\Ex\left[\norm{u_{n}(t\land \tau_N^n)
-u(t\land \tau_N^n)}^2_{H^\sigma}\Big|\F_0\right]\\
\leq \ &   
\norm{u_{0,n}-u_{0}}^2_{H^\sigma}
+\Ex\left[\int_0^{t\land \tau_N^n}\big(\Xi+CN+K(t',N)\big)
\norm{u_{n}(t')-u(t')}^2_{H^\sigma}\d t'\bigg|\F_0\right].
\end{align*}
Therefore, 
we can find a function $a(t,N)$, 
non-decreasing in $N$ and locally 
integrable in $t$ such that for  $t>0$, $n,\,N\ge1$,
\begin{align*}
&\Ex\left[
1\land \norm{u_{n}(t)-u(t)}^2_{H^\sigma}\Big|\F_0\right]\\
=\ &\Ex\left[
\left(1\land \norm{u_{n}(t)-u(t)}^2_{H^\sigma}\right)
\left(\mathbf{1}_{\{\tau_N^n<t\}}+\mathbf{1}_{\{\tau_N^n\ge t\}}\right)
\Big|\F_0\right]\\
\leq\ & \frac{\Big(V(\|u_{0,n}\|^2_{H^s})+V(\|u_0\|^2_{H^s})\Big){\rm e}^{\int_0^tg(t')\d t'}}{V(N^2)}
+\Big[1\land a(t,N)\|u_{0,n}-u_0\|^2_{H^\sigma}\Big].
\end{align*}
Therefore, we have that for all $t>0$ and $n\ge1$,
\begin{align*}
&\Ex\left[1\land \norm{u_{n}(t)-u(t)}^2_{H^\sigma}\Big|\F_0\right] \\
\leq\ &   \inf_{N\ge1}\left\{\frac{\Big(V(\|u_{0,n}\|^2_{H^s})+V(\|u_0\|^2_{H^s})\Big){\rm e}^{\int_0^tg(t')\d t'}}{V(N^2)}
+\Big(1\land a(t,N)\|u_{0,n}-u_0\|^2_{H^\sigma}\Big)\right\}.
\end{align*}
As a result, 
$$
\lim_{n\to\infty}\Ex
\left[
1\land \norm{u_{n}(t)-u(t)}^2_{H^\sigma}\Big|\F_0
\right]=0\ \  \pas,
$$
which implies \eqref{H-sigma stabilty}.
\end{proof}

\subsection{An evolution system of measures}\label{Section:T2:invariant measures}
In this section, we focus on the Cauchy problem \eqref{Cauchy problem-Ito} defined on the torus $\K=\T$.   
For any starting time $t'$ and target time $t$ such that $t' \le t$, we recall the Markov propagator $\mathscr{P}_{t',t}$ given by \eqref{eq:markovsg}, and the space $X_{s,\sigma}$ defined in \eqref{X-s1s2}. For $R \geq 1$, we define the closed ball
\begin{align}\label{eq:compact B-s-R}
	B^s_R \triangleq \big\{f \in H^s : \|f\|_{H^s} \le R\big\}.
\end{align}
From  Remark \ref{Remark: t0 in R} and \ref{T2:stability} in Theorem \ref{Thm: global 2}, one can show that (cf. \cite[Lemma 7.2]{Karlsen-Tang-Wang-2026-arXiv}) for any $t\ge t'$, $R>0$, and $\varphi \in C_B(X_{s,\sigma})$, the restricted mapping 
\begin{equation}
	\mathscr{P}_{t',t}\varphi\big|_{B^s_R}\  \text{is continuous with respect to the}\ H^{\sigma}\text{-topology}.\label{eq: restricted+mismatched feller}
\end{equation}

By initiating the system from the fixed deterministic state $u_0$ at varying past times $t'$, we introduce the approximate probability measures in $\mathbf{P}(H^{\sigma})$ via time-averaging:
\begin{align}\label{eq:approx_measures}
	\nu_{T,-n} \triangleq  \intbar_{-T}^{-n} \mathscr{P}^*_{t',-n} \delta_{u_0}\d t',\quad T > n \ge 0.
\end{align}

To establish \ref{T2:invariant measure} in Theorem \ref{Thm: global 2}, we use Lemma \ref{convergence of measures lemma} to show that the set $\{\nu_{T,-n}\}_{T > n} \subset \mathbf{P}(H^{\sigma})$ is tight. Taking the limit, we identify a candidate for an evolution system of measures in $\mathbf{P}(H^\sigma)$. Subsequently, we prove that this candidate indeed resides in the smaller space $\mathbf{P}(X_{s,\sigma})$, and that the weak convergence can be strengthened to the space $\mathbf{P}(X_{s,\sigma})$. These properties allow us to verify property \eqref{eq:limitproperty}.

\begin{proof}[Proof for \ref{T2:invariant measure} in Theorem \ref{Thm: global 2}]
	As outlined above, the proof proceeds in four steps.
	
	\textbf{Step 1: Tightness in $\mathbf{P}(H^{\sigma})$.} 
	Recall the solution operator $\mathscr{S}_{t',t}$ defined in \eqref{Solution operator}.  To establish the existence of an evolution system of measures, it suffices to consider a deterministic initial state. Adapting the proof of the Lyapunov estimate \eqref{Hs bound V tau-N} for any starting time $t_0=t'$ (see Remark \ref{Remark: t0 in R}), the expected value of the Lyapunov functional of the solution starting from the deterministic state $u_0$ at any past time $t' \le -n$ evaluated at time $-n$ satisfies:
	\begin{equation*}
		\mathbb{E}\left[V(\|\mathscr{S}_{t',-n} u_0\|_{H^s}^2)\right]  
		\leq V(\| u_0\|^2_{H^s}) \mathrm{e}^{\int_{t'}^{-n} g(r) \, \mathrm{d} r}.
	\end{equation*}

	Noting that $g$ in Hypothesis \ref{H-h-large} satisfies $g \in L^1(\mathbb{R};(0,\infty))$, we derive that
	\begin{equation}\label{Hs bound V t0}
		\mathbb{E}\left[V(\|\mathscr{S}_{t',-n} u_0\|_{H^s}^2)\right] \leq V(\| u_0\|^2_{H^s}) \mathrm{e}^{\|g\|_{L^1(\mathbb{R})}}, \quad \text{for all } t' \le -n.
	\end{equation}

Recall that $V' > 0$ (see \eqref{SCRV}), so $V$ is invertible. 
By the compact embedding $H^s \hookrightarrow\hookrightarrow H^\sigma$ on the torus $\T$, $B^s_R$ defined in \eqref{eq:compact B-s-R} is a compact subset of $H^{\sigma}$. This fact, combined with \eqref{Hs bound V t0}, implies that for any fixed $n$, the family $\{\nu_{T,-n}\}_{T> n}$ is tight in $\mathbf{P}(H^{\sigma})$.
Specifically, for any $\eta>0$, we can choose a sufficiently large radius $R = R_{\eta}> 0$ such that $V(R_{\eta}^2) >\eta^{-1} V(\|u_0\|^2_{H^s}) \exp{\{\|g\|_{L^1(\mathbb{R})}\}}$. By the definition of $\nu_{T,-n}$ in \eqref{eq:approx_measures}, and applying Chebyshev's inequality along with \eqref{eq:approx_measures} and \eqref{Hs bound V t0}, we obtain that for $u$ starting from $u_0$,	
\begin{align}
	\nu_{T,-n}(H^{\sigma}\setminus B^s_{R_{\eta}}) 
 \leq \frac{1}{V(R_{\eta}^2)} \int_{H^{\sigma}} V(\|v\|_{H^s}^2) \nu_{T,-n}(\mathrm{d} v)
	 = \frac{1}{V(R_{\eta}^2)} \intbar_{-T}^{-n} \mathbb{E}\left[V(\|\mathscr{S}_{t',-n}u_0\|_{H^s}^2)\right]\, \mathrm{d}t'  
< \eta, \label{tight of nu-(T,-n)}
\end{align}
which verifies the tightness of the family $\{\nu_{T,-n}\}_{T> n}$ in $\mathbf{P}(H^{\sigma})$.
As a result, we can infer from Lemma \ref{convergence of measures lemma} and Cantor's diagonal argument that there exists a common sequence of times $T_N \to \infty$ such that for each $n \in \mathbb{N}_0$, there is a probability measure $\nu_{-n} \in \mathbf{P}(H^{\sigma})$ satisfying $\nu_{T_N, -n}\to \nu_{-n}$ weakly in $\mathbf{P}(H^{\sigma})$, i.e.,
\begin{equation}\label{nu convergence Hsigma}
	\nu_{T_N, -n}\to \nu_{-n}\ \text{weakly in}\ \mathbf{P}(H^{\sigma}).
\end{equation}

\textbf{Step 2: Improving the convergence in $\mathbf{P}(X_{s,\sigma})$.} First, we assert that $\nu_{-n} \in \mathbf{P}(X_{s,\sigma})$. To demonstrate this, it suffices to establish that $\nu_{-n}$ is concentrated on the set $H^s$.  Using Lemma \ref{convergence of measures lemma} and the analysis in \eqref{tight of nu-(T,-n)}, we observe that
\begin{align*}
\nu_{-n}(H^{\sigma}\backslash B^s_R) \le 
\liminf_{N \uparrow \infty} \nu_{T_N, -n}(H^{\sigma}\backslash B^s_R) 
= o_{R \uparrow \infty}(1).
\end{align*}
Consequently,  from $\lim_{R \uparrow \infty} V(R) = \infty$, we conclude that
\begin{align*}
\lim_{R \uparrow \infty} \nu_{-n}(H^{\sigma}\setminus B^s_R)
= \nu_{-n}\bigg(\lim_{k \uparrow \infty} \bigcap_{R = 1}^k \big(H^{\sigma}\setminus B^s_R\big)\bigg)  = \nu_{-n}(H^\sigma \setminus H^s) = 0,
\end{align*}
indicating that $\nu_{-n}$ has zero mass outside the set $H^s$. Since the Borel $\sigma$-algebras on $X_{s,\sigma}$ and $H^s$ coincide, $\nu_{-n}$ is a well-defined probability measure on $H^s$.

It is clear that $\nu_{T_N, -n}\in \mathbf{P}(X_{s,\sigma})$. Since $X_{s,\sigma}\subset H^{\sigma}$ is dense, any function $f\in C_{B,U}(X_{s,\sigma})$ can be uniquely extended to a function $\widetilde{f}\in C_{B,U}(H^{\sigma})$. We can infer from this extension,  the fact $\nu_{-n}(H^\sigma\setminus H^s)=\nu_{T_N, -n}(H^\sigma\setminus H^s)=0$, and Lemma \ref{convergence of measures lemma}, that the convergence \eqref{nu convergence Hsigma} can be improved to
$
\nu_{T_N, -n}\to \nu_{-n}\ \text{weakly in}\ \mathbf{P}(X_{s,\sigma}),
$
that is,
\begin{equation}
\lim_{N\to\infty}\int_{H^s}f(v)\nu_{T_N, -n}({\rm d}v)=\int_{H^s}f(v)\nu_{-n}({\rm d}v),\quad f\in C_B(X_{s,\sigma}).\label{nu convergence X}
\end{equation}

\textbf{Step 3: Candidate evolution system of measures.}
Now that we have established $\nu_{-n} \in \mathbf{P}(X_{s,\sigma})$ and $\nu_{T_N, -n} \to \nu_{-n}$ weakly in $\mathbf{P}(X_{s,\sigma})$ (from the above step), we can define the candidate evolution system of measures. For any $n \ge 0$ such that $-n < t$, we define
\begin{align*}
	\mu_t^{(n)} \triangleq \mathscr{P}_{-n,t}^* \nu_{-n}, \quad t > -n.
\end{align*}
Since $\nu_{-n}$ is a well-defined probability measure on $H^s$, so is $\mu_t^{(n)}$. However, one must prove that $\mu_t^{(n)}$ is well-defined and independent of the starting time index $n$. To this end, we will first prove that
\begin{equation}\label{eq:rigorous_limit}
	\int_{H^s} \varphi(u)\mu^{(n)}_t(\mathrm{d} u) = \lim_{N \to \infty} \int_{H^s} [\mathscr{P}_{-n,t}\varphi](u) \nu_{T_N,-n}(\mathrm{d}u),\quad \varphi \in C_B(X_{s,\sigma}).
\end{equation}
Let $\psi(u) \triangleq [\mathscr{P}_{-n,t}\varphi](u)$. Due to \eqref{eq: restricted+mismatched feller}, $\psi$ is continuous with respect to the $H^\sigma$-topology only when restricted to bounded norm balls $B^s_{R}$ in $H^s$ (see \eqref{eq:compact B-s-R}). Since $\psi$ may not globally belong to $C_B(X_{s,\sigma})$, we cannot directly pass the weak convergence limit $N \to \infty$ into the integral $\int_{H^s} \psi(u) \nu_{T_N,-n}(\mathrm{d}u)$.

To justify \eqref{eq:rigorous_limit}, we borrow the idea in  \cite[Section 6.2]{Karlsen-Tang-Wang-2026-arXiv} on the time-homogeneous case. Fix an arbitrary $\eta > 0$. As in the analysis of \eqref{tight of nu-(T,-n)}, since $B^s_{R}$  is a compact set in $X_{s,\sigma}$ for any $R>0$, there exists a sufficiently large radius $R_0 > 0$ such that 
\begin{equation}\label{eq:tail_control_time}
	\sup_{N \ge 1} \nu_{T_N,-n}(H^s \setminus B^s_{R_0}) < \eta.
\end{equation}
The lower semicontinuity of the norm $\|\cdot\|_{H^s}$ with respect to the $H^\sigma$-topology ensures that $B^s_{R_0}$ is a closed subset of the metrizable space $X_{s,\sigma}$. By the Tietze extension theorem, there exists a globally continuous function $\Phi_{R_0} \in C_B(X_{s,\sigma})$ such that $\Phi_{R_0} \equiv \psi$ on $B^s_{R_0}$ and $\sup_{u \in H^s} |\Phi_{R_0}(u)| \le \sup_{u \in B^s_{R_0}} |\psi(u)| \le \|\varphi\|_{L^\infty}$.

We decompose the integration error as follows:
\begin{align*}
	&\left| \int_{H^s} \psi(u) \, \nu_{T_N,-n}(\mathrm{d}u)- \int_{H^s} \psi(u) \,  \nu_{-n}(\mathrm{d}u)\right| \notag\\
	\le\ &  \int_{H^s \setminus B^s_{R_0}} |\psi- \Phi_{R_0}|(u) \, \nu_{T_N,-n}(\mathrm{d}u)
+ \left| \int_{H^s} \Phi_{R_0}(u) \, \nu_{T_N,-n}(\mathrm{d}u) - \int_{H^s} \Phi_{R_0}(u) \, \nu_{-n}(\mathrm{d}u) \right|  
+\int_{H^s \setminus B^s_{R_0}} |\Phi_{R_0}- \psi|(u) \, \nu_{-n}(\mathrm{d}u).
\end{align*}
The first term is bounded by $2\|\varphi\|_{L^\infty}\eta $ due to \eqref{eq:tail_control_time}. For the third term, since the complement $H^s \setminus B^s_{R_0}$ is an open set in $X_{s,\sigma}$, 
Lemma \ref{convergence of measures lemma} implies $\nu_{-n}(H^s \setminus B^s_{R_0}) \le \liminf_{N \to \infty} \nu_{T_N,-n}(H^s \setminus B^s_{R_0}) \leq \eta $, bounding this term also by $2\|\varphi\|_{L^\infty} \eta $. For the middle term, the global continuity of $\Phi_{R_0} \in C_B(X_{s,\sigma})$ and the weak convergence $\nu_{T_N,-n} \rightharpoonup \nu_{-n}$ in $\mathbf{P}(X_{s,\sigma})$ ensure that it vanishes as $N \to \infty$. 

Taking the limit superior as $N \to \infty$, we obtain:
\begin{equation*}
	\limsup_{N \to \infty} \left| \int_{H^s} [\mathscr{P}_{-n,t}\varphi](u) \, \nu_{T_N,-n}(\mathrm{d}u) - \int_{H^s} [\mathscr{P}_{-n,t}\varphi](u) \, \nu_{-n}(\mathrm{d}u) \right| \le 4\|\varphi\|_{L^\infty} \eta .
\end{equation*}
Since $\eta> 0$ is arbitrary, we obtain \eqref{eq:rigorous_limit}.

Now, substituting the definition of $\nu_{T_N,-n}$ and employing the property $\mathscr{P}_{t',-n} \circ \mathscr{P}_{-n,t} = \mathscr{P}_{t',t}$, we deduce:
\begin{align*}
	\int_{H^s} \varphi(u)\mu^{(n)}_t(\mathrm{d} u) = \lim_{N \to \infty} \frac{1}{T_N - n} \int_{-T_N}^{-n} \int_{H^s} [\mathscr{P}_{t',-n} \circ \mathscr{P}_{-n,t}\varphi](u) \delta_{u_0}(\mathrm{d}u) \, \mathrm{d} t' 
	= \lim_{N \to \infty} \frac{1}{T_N - n} \int_{-T_N}^{-n} \mathscr{P}_{t',t} \varphi(u_0) \, \mathrm{d} t'.
\end{align*}

To verify that the evolution system is indeed independent of the starting backward time $n$, let us fix another integer $m > n$. A parallel argument shows that:
\begin{equation*}
	\int_{H^s} \varphi(u)\mu^{(m)}_t(\mathrm{d} u) = \lim_{N \to \infty} \frac{1}{T_N - m} \int_{-T_N}^{-m} \mathscr{P}_{t',t} \varphi(u_0) \, \mathrm{d} t'.
\end{equation*}
Since the integrand $\mathscr{P}_{t',t} \varphi(u_0)$ is globally bounded by $\|\varphi\|_{L^\infty}$, we can bound the difference between the time averages for $T_N > m$:
\begin{align*}
	&\left| \frac{1}{T_N - n} \int_{-T_N}^{-n} \mathscr{P}_{t',t} \varphi(u_0) \, \mathrm{d} t' - \frac{1}{T_N - m} \int_{-T_N}^{-m} \mathscr{P}_{t',t} \varphi(u_0) \, \mathrm{d} t' \right| \notag \\
	=\ & \left| \frac{1}{T_N - n} \int_{-m}^{-n} \mathscr{P}_{t',t} \varphi(u_0) \, \mathrm{d} t' + \left( \frac{1}{T_N - n} - \frac{1}{T_N - m} \right) \int_{-T_N}^{-m} \mathscr{P}_{t',t} \varphi(u_0) \, \mathrm{d} t' \right| \notag \\
	\le\ & \frac{m-n}{T_N - n} \|\varphi\|_{L^\infty} + \frac{m-n}{(T_N - n)(T_N - m)} (T_N - m) \|\varphi\|_{L^\infty} = \frac{2(m-n)}{T_N - n} \|\varphi\|_{L^\infty}.
\end{align*}
As $N \to \infty$ (and thus $T_N \to \infty$), this difference tends to $0$. This proves that 
\begin{equation}
	\int_{H^s} \varphi(u)\mu^{(n)}_t(\mathrm{d} u) = \int_{H^s} \varphi(u)\mu^{(m)}_t(\mathrm{d} u),\label{eq: mu-t n-independence}
\end{equation}
confirming that $\mu_t^{(n)}$ is well-defined and independent of $n$. Therefore, from now on we simply denote it as $\mu_t$.

\textbf{Step 4: Verification of \eqref{eq:limitproperty} for $\mu_t$.}  
To prove \eqref{eq:limitproperty}, we first
show that, for $-\infty < t' \leq t < \infty$, one may select $n>-t'\ge -t$ such that   for $\varphi\in C_B(X_{s,\sigma})$,
\begin{align} 
\int_{H^{s}}[\mathscr{P}_{t',t}\varphi](v)\mu_{t'}({\rm d}v)
= \int_{H^{s}}[\mathscr{P}_{-n,t'}\circ \mathscr{P}_{t',t}\varphi](v)\nu_{-n}({\rm d}v) 
=\int_{H^{s}}[\mathscr{P}_{-n,t}\varphi](v)\nu_{-n}({\rm d}v)
=\int_{H^{s}}\varphi(v)\mu_{t}({\rm d}v).\label{IM CB(X)}
\end{align}
We note that, for each $\phi\in C_B(H^s)$, there is a sequence $\{\phi_k\}_{k\ge1}\subset C_B(X_{s,\sigma})$ such that 
$$\sup_{k\ge1}\sup_{v\in H^s}|\phi_k(v)|<\infty,\quad \phi_k(v)\to\phi(v)\ \text{as}\ k\to\infty,\quad v\in H^{s}.$$
Indeed,  $\phi_k(\cdot)\triangleq \phi(J_k\cdot)$ with $J_k$ given by \eqref{Define Jn} verifies the above requirement.  From this and \eqref{IM CB(X)}, we obtain
\begin{align*}
	\int_{H^{s}}[\mathscr{P}_{t',t}\phi](v)\mu_{t'}({\rm d}v)
	=   \lim_{k\to\infty}\int_{H^{s}}[\mathscr{P}_{t',t}\phi_k](v)\mu_{t'}({\rm d}v) 
	=  \lim_{k\to\infty}\int_{H^s}\phi_k(v)\mu_{t}({\rm d}v)
	= \int_{H^s}\phi(v)\mu_{t}({\rm d}v),\quad \phi\in C_B(H^{s}).
\end{align*}
Hence, $\{\mu_t\}_{t \in \mathbb{R}}$ constitutes an evolution system of measures in the sense of Definition \ref{def:im_family}.
\end{proof}

\section*{Acknowledgement}

D.~A-O is supported by RYC2023-045563-I, funded by MCIU/AEI/10.13039/
501100011033 and FSE+, and partially by PID2023-148028NB-I00 funded by
MICIU/AEI/10.13039/501100011033 and also FEDER, UE.
P. P. was supported in part by the Research Council of Norway through the project INICE (301538) and by the University of Nottingham Ningbo China via FoSE Seed Funding for new academics.
H.~T.  is supported by the National Natural Science Foundation of China under projects  24AAA00530 and 12531007. 
H.~T. is also profoundly grateful to Professor Feng-Yu Wang for many inspiring discussions.

\setlength{\bibsep}{0.43ex}

\end{document}